\documentclass[12pt]{amsart}
\usepackage[utf8]{inputenc}

\title{The Nash-Williams orientation theorem for graphs with countably many ends}
\author[A.~Assem]{Amena Assem}
\address{Department of Combinatorics and Optimization, University of Waterloo, ON, Canada}
\email{a36mahmo@uwaterloo.ca}
\author[M.Koloschin]{Marcel Koloschin}
\author[M.~Pitz]{Max Pitz}
\address{Universit\"at Hamburg, Department of Mathematics, Bundesstrasse 55 (Geomatikum), 20146 Hamburg, Germany}
\email{\{marcel.koloschin, max.pitz\}@uni-hamburg.de}

\date{}



\usepackage{amsmath,amssymb,amsthm}  
\usepackage{mathtools}
\usepackage{tikz}
\usetikzlibrary{decorations.pathreplacing,decorations.markings}
\usepackage{comment}
\usepackage{enumitem}

\tikzset{
  on each segment/.style={
    decorate,
    decoration={
      show path construction,
      moveto code={},
      lineto code={
        \path [#1]
        (\tikzinputsegmentfirst) -- (\tikzinputsegmentlast);
      },
      curveto code={
        \path [#1] (\tikzinputsegmentfirst)
        .. controls
        (\tikzinputsegmentsupporta) and (\tikzinputsegmentsupportb)
        ..
        (\tikzinputsegmentlast);
      },
      closepath code={
        \path [#1]
        (\tikzinputsegmentfirst) -- (\tikzinputsegmentlast);
      },
    },
  },
  mid arrow/.style={postaction={decorate,decoration={
        markings,
        mark=at position .5 with {\arrow[#1]{stealth}}
      }}},
}

\usetikzlibrary{arrows}
\tikzset{
  arrow/.pic={\path[tips,every arrow/.try,->,>=#1] (0,0) -- +(.1pt,0);},
  pics/arrow/.default={triangle 90}
}

\usepackage{xcolor} 	
\usepackage[unicode]{hyperref}
\hypersetup{
	colorlinks,
    linkcolor={red!60!black},
    citecolor={green!60!black},
    urlcolor={blue!60!black},
}
\usepackage[abbrev, msc-links]{amsrefs} 

\usepackage[utf8]{inputenc}
\usepackage[T1]{fontenc}
\usepackage{lmodern}
\usepackage[babel]{microtype}
\usepackage[english]{babel}

\usepackage{geometry}
\usepackage{enumitem}

\let\polishlcross=\l
\def\l{\ifmmode\ell\else\polishlcross\fi}

\let\emptyset=\varnothing

\let\theta=\vartheta
\let\rho=\varrho
\let\phi=\varphi

\newcommand{\script}{\mathcal}
\newcommand{\parentheses}[1]{{\left( {#1} \right)}}

\newcommand{\p}{\parentheses}

\newcommand{\Set}[1]{{\left\lbrace {#1} \right\rbrace}}

\def\set#1:#2{\Set{{#1} \colon {#2}}}

\DeclareFontFamily{U}  {MnSymbolC}{}
\DeclareSymbolFont{MnSyC}         {U}  {MnSymbolC}{m}{n}
\DeclareFontShape{U}{MnSymbolC}{m}{n}{
    <-6>  MnSymbolC5
   <6-7>  MnSymbolC6
   <7-8>  MnSymbolC7
   <8-9>  MnSymbolC8
   <9-10> MnSymbolC9
  <10-12> MnSymbolC10
  <12->   MnSymbolC12}{}
\DeclareMathSymbol{\powerset}{\mathord}{MnSyC}{180}

\theoremstyle{plain}
\newtheorem{thm}{Theorem}[section]

\newtheorem{prop}[thm]{Proposition}
\newtheorem{clm}[thm]{Claim}

\newtheorem{lemma}[thm]{Lemma}

\newtheorem{obs}[thm]{Observation}

\theoremstyle{definition}

\begin{document}
\begin{abstract}
Nash-Williams proved in 1960 that a finite graph admits a $k$-arc-connected orientation if and only if it is $2k$-edge-connected, and conjectured that the same result should hold for all infinite graphs, too. 

Progress on Nash-Williams's problem was made by C.~Thomassen, who proved in 2016 that all $8k$-edge-connected infinite graphs admit a $k$-arc connected orientation, and by the first author, who recently showed that edge-connectivity of $4k$ suffices for locally-finite, 1-ended graphs. 

In the present article, we establish the optimal bound $2k$ in Nash-Williams's conjecture for all locally finite graphs with countably many ends.
\end{abstract}

\maketitle

\section{Introduction}

All graphs in this paper may have parallel edges but no loops. A directed graph is \emph{$k$-arc-connected} if from any vertex $v$ to any other vertex $w$ of the graph there exist $k$ arc-disjoint forwards directed paths. Clearly, the underlying undirected graph of a $k$-arc-connected graph must be $2k$-edge-connected. 
The classic \emph{orientation theorem} of Nash-Williams from 1960 asserts that for finite graphs, also the following converse is true.

\begin{thm}[Nash-Williams \cite{nash1960orientations}]
\label{thm_nash1960orientations}
Every finite $2k$-edge-connected graph has a $k$-arc-connected orientation.
\end{thm}

In the same paper, Nash-Williams claimed that his result also  holds for infinite graphs -- but 10 years later, Nash-Williams retracted his claim in \cite[\S~8]{nash1969well}.
Despite significant effort, it has remained open ever since whether the orientation theorem holds for infinite graphs as well. 

So far, for arbitrary infinite graphs, only the case $k=1$ was known, proved by Egyed by a Zorn's lemma argument already in 1941  \cite{egyed1941ueber}.
To appreciate the difficulty of the general case, note that a priori it is not even clear whether \emph{any} sufficiently large {(finite)} edge-connectivity implies the existence of a $k$-arc-connected orientation. This is different for finite graphs, where a simple argument shows that every $4k$-edge-connected graph has a $k$-arc-connected orientation: By the Nash-Williams/Tutte tree packing theorem \cite[Corollary~2.4.2]{Bible}, any such graph has $2k$ edge-disjoint spanning trees, so after fixing a common root, we may simply orient half of the trees away from and the other half towards the root. This approach, however, is blocked for infinite graphs: there exist locally finite graphs of arbitrarily large finite (edge-)connectivity that do not even posses two edge-disjoint spanning trees \cite{aharoni1989infinite}.

Motivated by the above considerations, Thomassen has asked in 1985 whether there is a function $f \colon \mathbb{N} \to \mathbb{N}$ such that any $f(k)$-edge-connected graph has a $k$-arc-connected orientation \cite{thomassen1989configurations}. This conjecture has been featured again in \cite[Conjecture~8]{barat2010his}, where also a topological variation of the problem was suggested by allowing directed topological arcs in the space consisting of $G$ together with all its ends; this topological version has been recently solved 
by Jannasch \cite{jannasch2019topologische}.

More than 50 years after Nash-Williams finite orientation theorem and about 30 years after posing his own conjecture, C.~Thomassen achieved a breakthrough towards the orientation theorem by {first reducing the problem to the locally finite case, and then by} proving  that every locally finite, $8k$-edge-connected graph has a $k$-arc-connected orientation \cite{thomassen2016orientations}, giving $f(k) \leq 8k$. 
Recently, A.~Assem further improved this bound to $4k$ in the case of one-ended, locally finite graphs \cite{assem2023towards}. Their proof strategies employ successively more refined variants of Mader's \emph{lifting theorem} from 1978 \cite{mader1978reduction}, 1992 \cite{frank} and 2016 \cite{ok2016liftings} and 2022 \cite{assem2022analysis}, results that were certainly not available to Nash-Williams in 1960.

In the present note (see Section~\ref{sec_orientation}), we further improve on Thomassen's and Assem's arguments in order to get the $4k$ bound for all graphs, and the best possible bound of $2k$ for all locally finite graphs with countably many ends.

\section{Boundary-linked decompositions}

Our first ingredient 
is a  modification of a decomposition result due to Thomassen, developed for his proof that $8k$-edge connectivity yields a $k$-arc-connected orientation.

Recall that an \emph{end} $\omega$ of a graph $G$ is an equivalence class of rays, where two rays of $G$ are \emph{equivalent} if there are infinitely many vertex-disjoint paths between them in~$G$. 
Now let $G=(V,E)$ be a locally finite connected graph. The \emph{boundary} of a set of vertices $B$ is the collection of edges in $G$ with one endvertex in $B$ and the other one outside of $B$. A set of vertices $B \subset V$ is called \emph{boundary-linked} if the induced subgraph $G[B]$ together with its boundary has a collection of pairwise edge-disjoint equivalent rays $R_1,R_2,\ldots $ such that each edge in the boundary is the first edge of one of the rays $R_i$. If we want to refer to which end $\omega$ those rays belong, we say that $B$ is \emph{$\omega$-boundary-linked}.

\begin{thm}[C.~Thomassen {\cite[Theorem~1]{thomassen2016orientations}}]
\label{thomass_finite partition}
Let $G$ be a locally finite, connected graph. Then for every finite set of vertices $A'$ in $G$, there is a finite set of vertices  $A \supseteq A'$ such that $V(G) \setminus A$ can be partitioned into finitely many boundary-linked sets each with finite boundary in $G$.
\end{thm}

Note, however, that the boundary linked sets in Thomassen's theorem do not generally coincide with the components of $G-A$, i.e.\ there might be edges between two boundary-linked sets. Our next result shows that for graphs with only countably many ends, we can obtain such a stronger partition:

\begin{thm}\label{finite partition}
Let $G$ be a locally finite, connected graph with at most countably many ends. Then for every finite set of vertices $A'$ in $G$, there is a finite set of vertices $A \supseteq A'$  such that all components of $G-A$ are boundary-linked. 
\end{thm}

Before we turn towards the proof, we remark that the assumption of $G$ having only countably many ends is necessary. Indeed, the following graph contains no finite set of vertices $A$ of size at least $2$ such that all components of $G-A$ are boundary linked:
\begin{figure}[h]
    \centering
  \begin{tikzpicture}[scale=0.8]

  \node[circle, draw, inner sep=0pt,minimum width=3pt, fill] (a) at (0,0) {};

    \node[circle, draw, inner sep=0pt,minimum width=3pt, fill] (a0) at (-1.7,2) {};
     \node[circle, draw, inner sep=0pt,minimum width=3pt, fill] (a1) at (1.7,2) {};

        \node[circle, draw, inner sep=0pt,minimum width=3pt, fill] (a00) at (-3,4) {};
     \node[circle, draw, inner sep=0pt,minimum width=3pt, fill] (a01) at (-1,4) {};

       \node[circle, draw, inner sep=0pt,minimum width=3pt, fill] (a10) at (1,4) {};
     \node[circle, draw, inner sep=0pt,minimum width=3pt, fill] (a11) at (3,4) {};

          \node[circle, draw, inner sep=0pt,minimum width=3pt, fill] (a000) at (-3.5,6) {};
     \node[circle, draw, inner sep=0pt,minimum width=3pt, fill] (a001) at (-2.5,6) {};

       \node[circle, draw, inner sep=0pt,minimum width=3pt, fill] (a010) at (-1.5,6) {};
     \node[circle, draw, inner sep=0pt,minimum width=3pt, fill] (a011) at (-0.5,6) {};

              \node[circle, draw, inner sep=0pt,minimum width=3pt, fill] (a100) at (0.5,6) {};
     \node[circle, draw, inner sep=0pt,minimum width=3pt, fill] (a101) at (1.5,6) {};

       \node[circle, draw, inner sep=0pt,minimum width=3pt, fill] (a110) at (2.5,6) {};
     \node[circle, draw, inner sep=0pt,minimum width=3pt, fill] (a111) at (3.5,6) {};

     \draw (a) -- (a0);
     \draw (a) -- (a1);

     \draw (a0) -- (a00);
     \draw (a0) -- (a01);
     \draw (a1) -- (a10);
     \draw (a1) -- (a11);

     \draw (a00) -- (a000);
     \draw (a00) -- (a001);
     \draw (a01) -- (a010);
     \draw (a01) -- (a011);

     \draw (a10) -- (a100);
     \draw (a10) -- (a101);
     \draw (a11) -- (a110);
     \draw (a11) -- (a111);

     \begin{scope}[shift={(0.4,0)}]
         
  \node[circle, draw, inner sep=0pt,minimum width=3pt, fill] (b) at (0,0) {};

    \node[circle, draw, inner sep=0pt,minimum width=3pt, fill] (b0) at (-1.7,2) {};
     \node[circle, draw, inner sep=0pt,minimum width=3pt, fill] (b1) at (1.7,2) {};

        \node[circle, draw, inner sep=0pt,minimum width=3pt, fill] (b00) at (-3,4) {};
     \node[circle, draw, inner sep=0pt,minimum width=3pt, fill] (b01) at (-1,4) {};

       \node[circle, draw, inner sep=0pt,minimum width=3pt, fill] (b10) at (1,4) {};
     \node[circle, draw, inner sep=0pt,minimum width=3pt, fill] (b11) at (3,4) {};

          \node[circle, draw, inner sep=0pt,minimum width=3pt, fill] (b000) at (-3.5,6) {};
     \node[circle, draw, inner sep=0pt,minimum width=3pt, fill] (b001) at (-2.5,6) {};

       \node[circle, draw, inner sep=0pt,minimum width=3pt, fill] (b010) at (-1.5,6) {};
     \node[circle, draw, inner sep=0pt,minimum width=3pt, fill] (b011) at (-0.5,6) {};

              \node[circle, draw, inner sep=0pt,minimum width=3pt, fill] (b100) at (0.5,6) {};
     \node[circle, draw, inner sep=0pt,minimum width=3pt, fill] (b101) at (1.5,6) {};

       \node[circle, draw, inner sep=0pt,minimum width=3pt, fill] (b110) at (2.5,6) {};
     \node[circle, draw, inner sep=0pt,minimum width=3pt, fill] (b111) at (3.5,6) {};

     \draw (b) -- (b0);
     \draw (b) -- (b1);

     \draw (b0) -- (b00);
     \draw (b0) -- (b01);
     \draw (b1) -- (b10);
     \draw (b1) -- (b11);

     \draw (b00) -- (b000);
     \draw (b00) -- (b001);
     \draw (b01) -- (b010);
     \draw (b01) -- (b011);

     \draw (b10) -- (b100);
     \draw (b10) -- (b101);
     \draw (b11) -- (b110);
     \draw (b11) -- (b111);
     \end{scope}

     \draw (a) -- (b);
     \draw (a0) -- (b0) -- (a1) -- (b1);
      \draw (a00) -- (b00) -- (a01) -- (b01);
       \draw (a10) -- (b10) -- (a11) -- (b11);

        \draw (a000) -- (b000) -- (a001) -- (b001);
       \draw (a010) -- (b010) -- (a011) -- (b011);

        \draw (a100) -- (b100) -- (a101) -- (b101);
       \draw (a110) -- (b110) -- (a111) -- (b111);

    \node (text) at (-2,7.2) {$\vdots$};
     \node (text) at (2,7.2) {$\vdots$};
   \end{tikzpicture}
    \caption{A tree-like graph witnessing the optimality of Theorem~\ref{finite partition}.}
    \label{fig:counterexample}
\end{figure}
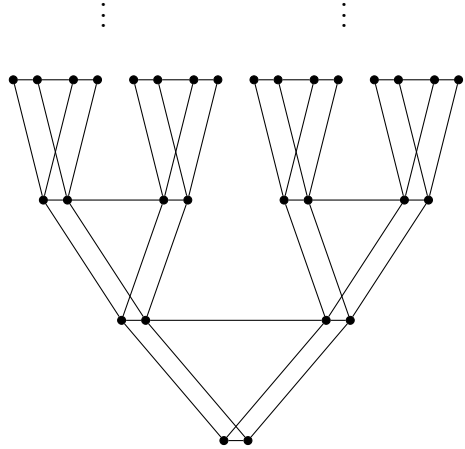

Formally, let $T$ be the rooted, infinite binary tree, whose vertices are given by the finite, binary sequences. Let $G_0 = T \times \{0\}$ and $G_1 = T \times \{1\}$ be two copies of the binary tree. Then $G$ is the graph consisting of $G_0$ and $G_1$ together with all edges of the form $(t,0)(t,1)$ for $t \in T$ as well as $(t^\frown 0,0)(t^\frown 1,1)$ for all $t \in T$, cf.~Figure~\ref{fig:counterexample}. Then every end of $G$ contains at most $3$ pairwise edge-disjoint rays; but for every finite set of vertices $A$ of size at least $2$, there will be a component of $G-A$ with boundary size at least $4$. We omit the details.

\begin{proof}[Proof of Theorem~\ref{finite partition}]
Let $\{\omega_1,\omega_2,...\}$ be an enumeration of the ends of $G$. Let $E_1$ be a cut of minimal size separating $A'$ from $\omega_1$. 
By minimality of $E_1$, the component $B_1$ of $G-E_1$ containing $\omega_1$ has boundary $E_1$, 
and is boundary-linked to $\omega_1$, see e.g.\ \cite[Lemma~10]{bruhn2007end}. Let $A_1 := N(B_1)$ consist of all vertices outside of $B_1$ with a neighbour in $B_1$.

Now assume that for $n\geq1$ we have already defined finite sets of vertices $A_1, \ldots, A_n$ and disjoint boundary-linked sets $B_1, \ldots, B_n$. Since $G$ is locally finite, the subgraph
$$G - \left(A' \cup \bigcup_{i=1}^{n} (A_i \cup B_i)\right)$$
has only finitely many components. If all of them are finite, then we add the remaining vertices to $A' \cup A_1 \cup \cdots \cup A_n$ to obtain our desired set of vertices $A$ with boundary-linked components $B_1,\ldots,B_n$. If not, then each infinite component contains a ray (again, because $G$ is locally finite). Let $\omega_{m}$ be the first end in our enumeration for which there is a ray left. 
As above, let $E_{n+1}$ be a cut of minimal size separating $A' \cup A_1 \cup \cdots \cup A_n$ from $\omega_m$. Then the component $B_{n+1}$ of $G-E_{n+1}$ containing $\omega_m$ has boundary $E_{n+1}$ and is boundary-linked to $\omega_m$. Letting $A_{n+1}$ be the endvertices of edges in $E_{n+1}$ outside of $B_{n+1}$ completes the recursive construction.

Suppose for a contradiction that this procedure does not terminate. In this case, the sets $B_1,B_2,...$ induce an open partition of the endspace of $G$. Since $G$ is locally finite, the endspace is compact \cite[Proposition~8.6.1]{Bible}, so there is a finite subcover. It follows that we already covered all ends after finitely many steps, a contradiction.
\end{proof}

\section{Immersions of finite  graphs  of prescribed connectivity}

If $G$ is a graph and $H$ is a graph with vertices $x_1,x_2,\ldots,x_n$, then an \emph{immersion} of $H$ in $G$ is a subgraph of $G$ consisting of $n$ distinguished vertices $y_1,y_2,\ldots,y_n$ and a collection of pairwise edge-disjoint paths in $G$ such that for each edge $x_ix_j$ in $H$ there is a corresponding path in the collection from $y_i$ to $y_j$. This immersion is said to be \emph{on} $\{y_1,\ldots,y_n\}$.
This subgraph of $G$ is a \emph{strong immersion} of $H$ if additionally, any such path joining $y_i,y_j$ avoids all other distinguished vertices $y_k$.

\begin{thm}[{C.~Thomassen \cite{thomassen2016orientations}*{Theorem~4}}]\label{immersion_thomassen}
Let $k$ be a positive integer, $G=(V,E)$ be a $2k$-edge-connected graph and $A \subseteq V$ be a finite set of vertices such that $V(G) \setminus A$ can be partitioned into finitely many boundary linked sets each with finite boundary.
Then $G$ contains a strong immersion of a finite Eulerian $k$-edge-connected graph on the vertex set $A$.
\end{thm}

Thomassen's proof proceeds as follows. First, he contracts each boundary linked set to a new dummy vertex, resulting in a finite $2k$-edge-connected graph $H'$. Then, Thomassen spends half the connectivity to move to a spanning Eulerian $k$-edge-connected subgraph $H$ of $H'$. To obtain the final Eulerian immersion on $A$, he successively replaces pairs of edges at dummy vertices by pairwise edge-disjoint paths inside the corresponding boundary linked sets such that $k$-edge-connectivity is preserved for vertices in $A$ throughout. For this last step, Thomassen proved a variant of Mader's lifting theorem, for Eulerian graphs, by studying what he called the \emph{bad graph}, sometimes known as the \emph{nonadmissibility graph}, which encodes the information of which pairs of edges incident with the dummy vertex are not replaceable by paths as described above (similar results were already proven by T.~Jordan in \cite[Theorem~3.2]{doi:10.1137/S0895480199364483}). 

To reach our target of optimizing the constants, however, we must make do without the Eulerian condition, which causes a number of technical challenges. Here is our eventual result that takes the role of Thomassen's Theorem~\ref{immersion_thomassen}. {Write $\lambda_G(x,y)$ for the local edge-connectivity between two vertices $x$ and $y$, i.e. the maximum number of edge-disjoint paths between two vertices $x$ and $y$ in $G$.}

\begin{thm}
\label{immersion}
Let $k\geq 2$ be an integer and $G$ be a $2k$-edge-connected, locally finite graph, and let $A$ be a set of vertices such that each component of $G-A$ is boundary linked. 

Then there is a set of vertices $X \subseteq V(G) \setminus A$ containing at most one vertex from each boundary-linked component such that $G$ contains an immersion of a finite, $3$-edge-connected graph $H$ on $A \cup X$ with the following properties:
\begin{itemize}
    \item [(i)] $d_H(x)=3$ for all $x \in X$,
    \item [(ii)] $\lambda_H(a,b) \geq 2k$ for all distinct $ a,b \in A$, and
    \item[(iii)] any edge in $G$ between two vertices in $A$ is also an edge in $H$.
\end{itemize}
\end{thm}

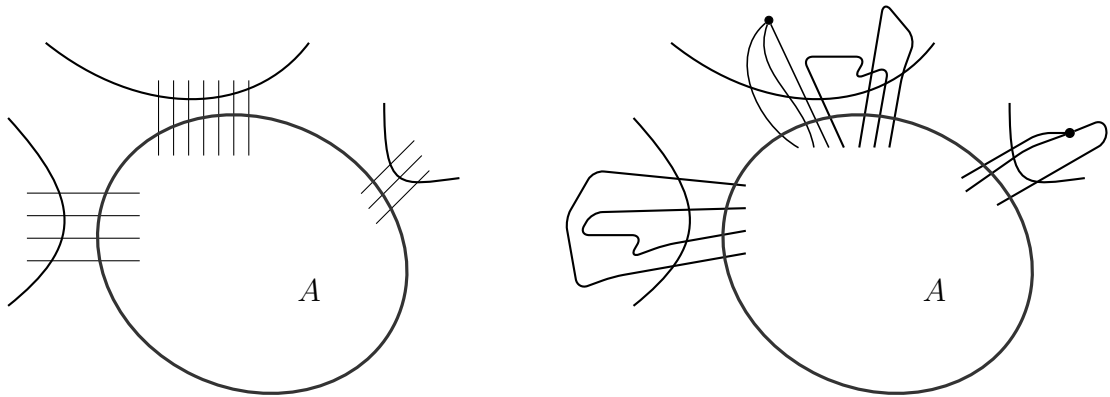
\begin{figure}[ht]
\centering
\begin{tikzpicture}
\tikzstyle{O}=[draw,very thick,black!30]
        \tikzstyle{dot}=[draw,circle,color=white,minimum size=0pt,inner sep=0pt]

\tikzstyle{white node}=[draw,circle,fill=white,minimum size=4pt, inner sep=0pt]
                            \tikzstyle{edge} = [draw,line width=0.5pt,-]
    \tikzstyle{green node}=[draw,circle,fill=green,minimum size=4pt, inner sep=0pt]
  \tikzstyle{edge} = [draw,line width=1pt,-]


\draw[thick] (-3.5,1.5) .. controls (-2.25,0.5) and (-0.75,0.5) .. (0,1.5);
\draw[thick] (-4,0.5) .. controls (-3,-0.6) and (-3,-1.2) .. (-4,-2);
\draw[thick] (1,0.7) .. controls (1,-0.5) and (1.2,-0.4) .. (2,-0.3);
        \draw (-3.75,-0.5) node[dot] () []{}
                    -- ++(180:-1.5cm) node[dot] () {};
        \draw (-3.75,-0.8) node[dot] () []{}
                    -- ++(180:-1.5cm) node[dot] () {};
        \draw (-3.75,-1.1) node[dot] () []{}
                    -- ++(180:-1.5cm) node[dot] () {};
         \draw (-3.75,-1.4) node[dot] () []{}
                    -- ++(180:-1.5cm) node[dot] () {};
                    
         \draw (-2,1) node[dot] () []{}
                    -- ++(90:-1cm) node[dot] () {};
         \draw (-1.8,1) node[dot] () []{}
                    -- ++(90:-1cm) node[dot] () {};
         \draw (-1.6,1) node[dot] () []{}
                    -- ++(90:-1cm) node[dot] () {};
         \draw (-1.4,1) node[dot] () []{}
                    -- ++(90:-1cm) node[dot] () {};
         \draw (-1.2,1) node[dot] () []{}
                    -- ++(90:-1cm) node[dot] () {};
         \draw (-1,1) node[dot] () []{}
                    -- ++(90:-1cm) node[dot] () {};
         \draw (-0.8,1) node[dot] () []{}
                    -- ++(90:-1cm) node[dot] () {};
                    
         \draw (1.5,0) node[dot] () []{}
                    -- ++(45:-1cm) node[dot] () {};
         \draw (1.4,0.2) node[dot] () []{}
                    -- ++(45:-1cm) node[dot] () {};
         \draw (1.6,-0.2) node[dot] () []{}
                    -- ++(45:-1cm) node[dot] () {};
          
 \draw[rotate around={-25:(0.9,2)},very thick,black!80] (0.8,-1.7) ellipse (60pt and 51pt);

        \draw (0,-1.5) node[dot] () [label=-90:$A$]{};

\end{tikzpicture}
\quad \quad \quad
\begin{tikzpicture}
\tikzstyle{O}=[draw,very thick,black!30]
        \tikzstyle{dot}=[draw,circle,color=white,minimum size=0pt,inner sep=0pt]
                            \tikzstyle{edge} = [draw,line width=0.5pt,-]
    \tikzstyle{green node}=[draw,circle,fill=black,minimum size=4pt,
                            inner sep=0pt]
  \tikzstyle{edge} = [draw,line width=0.6pt,-]

    \tikzstyle{black node}=[draw,circle,fill=black,minimum size=3pt,
                            inner sep=0pt]

\draw[thick] (-3.5,1.5) .. controls (-2.25,0.5) and (-0.75,0.5) .. (0,1.5);
\draw[thick] (-4,0.5) .. controls (-3,-0.6) and (-3,-1.2) .. (-4,-2);
\draw[thick] (1,0.7) .. controls (1,-0.5) and (1.2,-0.4) .. (2,-0.3);

                            \draw (-2.5,-1) node[dot] (x1) []{};
                            \draw (-2.5,-1.3) node[dot] (x2) []{};
                            \draw (-2.5,-0.7) node[dot] (x3) []{};
                            \draw (-2.5,-0.4) node[dot] (x4) []{};

   \draw[rounded corners, thick] (x2)  -- ++(190:1.8cm) -- ++(200:0.5cm) -- ++(100:0.95cm)-- ++(60:0.75cm) -- ++(x4);
   \draw[rounded corners, thick] (x1)  -- ++(190:1.1cm) -- ++(200:0.5cm) -- ++(60:0.35cm) -- ++(180:0.85cm) -- ++(60:0.35cm) -- ++(x3);
                 \draw (-2.2,1.8) node[black node] (y) []{};
                                    
                \draw (-0.6,0.1) node[dot] (y1) []{};
                \draw (-0.8,0.1) node[dot] (y2) []{};
                \draw (-1,0.1) node[dot] (y3) []{};
                \draw (-1.2,0.1) node[dot] (y4) []{};
                \draw (-1.4,0.1) node[dot] (y5) []{};
                \draw (-1.6,0.1) node[dot] (y6) []{};
                \draw (-1.8,0.1) node[dot] (y7) []{};
                                                        
        \path[edge] (y) to (y5);
        \path[edge] (y) to [out=250,in=100] (y6);
        \path[edge] (y) to [out=220,in=150] (y7);       

    \draw[rounded corners, thick] (y1)  -- ++(80:1cm) -- ++(70:0.5cm) -- ++(130:0.65cm) -- ++(y3);
   \draw[rounded corners, thick] (y2)  -- ++(80:1.1cm) -- ++(200:0.5cm) -- ++(60:0.35cm) -- ++(180:0.85cm) -- ++(y4);
                \draw[rounded corners, thick] (1.8,0.3) node[black node] (z) []{}
                     -- ++(180:0.4cm)  -- ++(30:-1.2cm) node[dot] () {};
                \draw[rounded corners, thick] (z)
                    -- ++(20:-0.6cm) -- ++(35:-1cm) node[dot] () {};
                \draw[rounded corners, thick] (z)
                     -- ++(22:0.5cm) -- ++(-82:0.3cm) -- ++(30:-1.7cm) node[dot] () {};

 \draw[rotate around={-25:(0.9,2)},very thick,black!80] (0.8,-1.7) ellipse (60pt and 51pt);
 
 \draw (0,-1.5) node[dot] () [label=-90:$A$]{};
 
 \end{tikzpicture}
 \caption{Connecting pairs of edges.}
\end{figure}

We picture this as follows: Starting from a family of boundary linked sets as illustrated on the left, we carefully connect pairs of edges by edge-disjoint paths in $G$ (each represented by an additional edge in $H$) until -- for boundary-linked components with odd boundary -- there are 3 edges left, which we connect via edge-disjoint paths to a vertex in $X$ resulting in a vertex of degree $3$ in $H$, as pictured on the right (in the immersed graph, $x$ might  have a larger, odd degree). 

In order to make this precise, we need to introduce the lifting operation pioneered by Lov{\'a}sz \cite{lovasz1976some}, Mader \cite{mader1978reduction} and Frank \cite{frank}. Let $G=(V+s,E)$ be a finite graph. \emph{Lifting} two distinct edges  $sx,sy$ incident with $s$ means deleting them and adding instead a new edge $xy$ (possibly parallel to existing edges between $x$ and $y$). If $x=y$, we delete the resulting loop.
%

A \textit{target connectivity} is a function $\tau$, whose values are non-negative integers, such that $\tau(x,y)\leq \lambda_{G}(x,y)$ for any two vertices $x$ and $y$ in $G$.
Suppose $G=(V+s,E)$ is a finite graph.  A pair of edges incident with $s$ is called \emph{$\tau$-admissible} if after lifting them, we obtain a graph $G'$ such that $\lambda_{G'}(x,y) \geq \tau(x,y)$ for all $x,y \in V$. 
{When the target connectivity $\tau$ is understood from the context, we simply write \emph{admissible}, and we say that an edge is \emph{liftable} with another edge if the pair is admissible.}
The corresponding \emph{lifting graph} $L(G,s,\tau)$ is the graph whose vertices are the edges incident with $s$, where two of them are adjacent if and only if they form a $\tau$-admissible pair. 

In \cites{thomassen2016orientations, doi:10.1137/S0895480199364483, ok2016liftings, assem2022analysis}, the lifting graph was studied for $k$-edge-connected graphs and constant target function $\tau \equiv k$. In the context of Theorem~\ref{immersion}, however, we have a prescribed vertex set $A$ and several (dummy) vertices for which, one after another, we want to lift edges in order to find the desired immersion. In this context, we are  interested only in preserving the edge-connectivity $2k$ for vertices in $A$.
But this requirement can easily be modelled in terms of a target function $\tau$. Indeed, define
$$\tau_A(x,y) = \begin{cases} 2k & \text{ if } x,y \in A, \\ 0 & \text{ otherwise.} \end{cases}$$

Now our proof relies on the following structural result on the lifting graph, which fine tunes the methods from \cite{assem2022analysis,ok2016liftings}.

\begin{thm} \label{main lifting theorem}
Let $G$ be a finite graph, {$k$ a positive}  integer, and $A$ a proper subset of $V(G)$ such that between any two vertices in $A$ there are $2k$ edge-disjoint paths in $G$. Assume that there is no edge with both endvertices outside of $A$, and let $s$ be a vertex not in $A$ with $\deg(s) > 3$. Then $L(G,s,\tau_A)$ is either
\begin{itemize}

\item a complete multipartite graph, or

\item {the union of an isolated vertex and a complete bipartite graph with equally sized partition classes. Hence, this latter case can only occur if $\deg(s)$ is odd.}
\end{itemize}

\end{thm}

We postpone the proof of Theorem~\ref{main lifting theorem} to the next section, so that we can immediately show how Theorem \ref{immersion} follows. When speaking of an admissible pair of edges, it is always understood that we mean admissible with respect to $\tau_A$. In the proof below, we will make use of the following well-known observation:

\begin{obs}
\label{obs_does not become feasible again}
 If, after lifting the admissible pair $\Set{e_1, e_2}$, the pair $\Set{e_3, e_4}$ is admissible, then ${e_3, e_4}$ is admissible in the original graph.
\end{obs}

\begin{proof}[Proof of Theorem~\ref{immersion}]
Given Theorem~\ref{main lifting theorem}, the proof now follows almost exactly along the lines of \cite[Lemma~3.1]{assem2023towards}. We give the full argument for the convenience of the reader. 

So let $G$ be a $2k$-edge-connected, locally finite graph, and let $A$ be a set of vertices such that each component of $G-A$ is boundary linked. Let $G_1$ be the finite graph obtained by contracting each of the boundary-linked components $B_1,\ldots,B_n$ into a single vertex $s_1,\ldots,s_n$ (keeping all parallel edges that arise). Then $G_1$ is $2k$-edge-connected. Since $k \geq 2$,  $\deg(s_i) > 3$ and $G_1$ has no bridge. Moreover, {$A \subseteq V(G_1)$} is such that every edge in $G_1$ is incident with a vertex from $A$.

We now lift edges at $s_1,\ldots,s_n$ in turn, preserving the property that $\lambda(x,y)\geq 2k$ for $x,y\in A$. Consider the lifting graph $L(G_1,s_1,\tau_A)$. Its vertices are the edges $e_1,e_2,...,e_{\deg(s_1)}$ incident with $s_1$.
We now define another graph {$M(\script{R})$}, the \emph{ray graph}, on this vertex set. We consider a set $\script{R}$ of edge-disjoint rays $R_1,R_2,...,R_{\deg(s_1)}$ in $G[B_1]$ starting with the edges $e_1,e_2,...,e_{\deg(s_1)}$ in the boundary of $B_1$, and add an edge between two vertices $e_i,e_j$ in {$M(\script{R})$} if in $G[B_1]$ there are
infinitely many pairwise disjoint paths joining $R_i,R_j$ having only their endpoints in
common with $\bigcup \script{R}$. Since these rays belong to the same end of $G$, it follows that {$M(\script{R})$} is
connected. 

Now if the two graphs $L(G_1,s_1,\tau_A)$ and {$M(\script{R})$} have a common edge joining $e_i,e_j$, say, then let $P$ be a path in $G[B_1]$ joining $R_i,R_j$ such that $P$ has only its endvertices in common with $\bigcup \script{R}$. Let $P_{i,j}$ be a path in $R_i \cup R_j \cup P$ starting and terminating
with $e_i,e_j$. Now delete the edges of $P_{i,j}$ from $G[B_1]$, lift $e_i,e_j$ in $G_1$ resulting in a new graph $G'_1$ and
define a new ray graph {$M(\mathcal{R}\setminus \{R_i,R_j\})$} 
and a new lifting graph $L(G'_1,s_1,\tau_A)$ (without the edges $e_i,e_j$). 
If the new ray graph and the new
lifting graph have a common edge, we repeat the above argument to find a new path in $G[B_1]-E(P_{i,j})$ and lift the corresponding edges in the new $G'_1$. Note that this new path may go through the rays disregarded from $\mathcal{R}$ ({but is edge-disjoint to the previously selected path as we have deleted the edges of that path}). Repeat this procedure for as long as possible. There are three possible outcomes:

\textbf{Case 1:} $\deg(s_1)$ is even. Then by Theorem~\ref{main lifting theorem}, the lifting graph is always a complete multipartite graph, so it has disconnected complement. Since $M$ is connected, the two graphs will always have a common edge, and so the above procedure will continue until all edges at $s_1$ have been lifted.

\textbf{Case 2:} $\deg(s_1)$ is odd, and the above procedure continues until $\deg(s_1) = 3$. Write $e_i,e_j$ and $e_k$ for the three remaining edges incident with $s_1$. Since the current ray graph {$M=M(\{R_i,R_j,R_k\})$} 
is connected, one of the edges, say $e_i$, is connected to the other two in $M$. Pick paths $P_j$ and $P_k$ in $G[B_1]$ joining $R_i,R_j$ and $R_i,R_k$. Inside $R_i \cup R_j \cup R_k \cup P_j \cup P_k$, we find a vertex $x_1$ together with three edge-disjoint paths from $x_1$ to $e_i,e_j,e_k$. 

\textbf{Case 3:} $\deg(s_1)$ is odd, and the above procedure terminates before $s_1$ reaches degree $3$. Then by Theorem~\ref{main lifting theorem}, the current lifting graph consists of an isolated vertex $e^*$ and a complete bipartite graph with partition classes say $e_1,\ldots,e_j$ and $e'_1,\ldots,e'_j$. {Let $M=M(\script{R'})$ be the current ray graph where $\script{R'}$ are the rays in $\script{R}$ starting in the edges of the current lifting graph.} Since $M$ is connected and {has no edge in common with the current lifting graph by Definition of Case~3,} 
$e^*$ has neighbours say $e_1$ and $e'_1$ in $M$. Let $R^*$, $R$ and $R'$ be the rays of $\mathcal{R}'$ that begin with $e^*$, $e_1$ and $e'_1$. 
Pick disjoint paths $P$ and $P'$ between $R,R^*$ and $R',R^*$ {that avoid all other rays in $\mathcal{R}'$}. Inside $R^* \cup R \cup R' \cup P \cup P'$, we find a vertex $x_1$ together with three edge-disjoint paths from $x_1$ to $e^*,e_1,e'_1$. Write $W_1$ for this path system; delete its edges.

Then define a new ray graph  {$M=M(\mathcal{R}'')$ on $\mathcal{R}'' =\mathcal{R}' \setminus \{R^*,R,R'\}$}, but keep the old lifting graph (including $e^*,e_1,e'_1$). Then the lifting graph has a common edge $e_i,e'_i$ with $M$ (because the vertices of $M$ induce a complete bipartite subgraph of the lifting graph, and $M$ is connected). As before, let $P$ be a path in $G[B_1]$ joining the corresponding rays $R_i,R'_i$ such that $P$ has only its endpoints in common with $\bigcup \script{R}''$. Let $P_i$ be a path in $R_i \cup R'_i \cup P$ starting and terminating
with $e_i,e'_i$. Now delete the edges of $P_i$ from $G[B_1]$, lift $e_i,e'_i$, and define a new ray graph $M$ (where we disregard $R_i, R'_i$) and a new lifting graph without the edges $e_i,e'_i$ (but still including $e^*,e_1,e'_1$). 

By Observation~\ref{obs_does not become feasible again}, $e^*$ is still isolated in this new lifting graph, and so by Theorem~\ref{main lifting theorem} the vertices of $M$ induce a complete bipartite subgraph of the lifting graph. Hence, they have a common edge, and we may proceed as before, until all edges but $e^*,e_1,e'_1$ have been lifted; but for these, we already reserved a vertex $x_1$ and a path system $W_1$. Note that at every step the new path connecting the two rays under consideration may go through the previously disregarded rays. Also the two rays being joined may be separated in the original ray graph only by the rays that were disregarded before connecting them.

\medskip
This concludes the lifting procedure at the first dummy vertex $s_1$.
Write $G_2$ for the resulting graph (after deleting $s_1$ in case all the edges incident with it are lifted). By construction, we have $\lambda_{G_2}(a,b) \geq 2k$ for all $a,b \in A$. Note also that we only lifted edges with one {endvertex} in $A$, so any edge with both {endvertices} in $A$ remains (this will ensure (iii) is satisfied).
Now consider the lifting graph $L(G_2,s_2,\tau_A)$, and proceed as above. Dealing with $s_3,\ldots,s_n$ in turn completes the construction. To see that the final graph is $3$-edge-connected, consider an edge-cut with one side not containing vertices from $A$. Then that side of the cut does not {span} any edges as there are no edges between the vertices $s_1,\ldots,s_n$ which resulted from contracting the components of $G-A$. So any vertex on that side of the cut contributes all of its incident edges to the cut, and the cut has size at least $3$.
\end{proof}

\section{Liftings}

In this section, we prove Theorem~\ref{main lifting theorem}. 
The proof is a routine adaption of the corresponding methods from \cite{assem2022analysis,ok2016liftings}. For the sake of completeness, we give the details in full. To make for less clustered reading, all results following Theorem~\ref{dangerous set theorem} below should be read under the assumptions of Theorem~\ref{main lifting theorem}.

Given a graph $G$ and sets of vertices $X,Y$ of $G$ we write $E(X,Y)$ for all edges with one {endvertex} in $X$ and the other in $Y$; and $\delta(X) = E(X,\overline{X})$ for the edge-boundary of $X$; here, $\overline{X}$ denotes the complement $V(G) \setminus X$ of $X$.

We will frequently use the following standard equation describing the intersection of two edge cuts $\delta(A_1)$ and $\delta(A_2)$ for any $A_1,A_2\subseteq V(G)$.

\begin{equation} \label{standard equation}
\begin{split}
2 \bigg[ &\big|\delta(A_1)\big| +\big |\delta(A_2)\big|- \Big (\big|E(A_1 \cap A_2, \overline{A_1 \cup A_2})\big|+\big|E(A_2\setminus A_1, A_1\setminus A_2)\big|\Big) \bigg] 
\\
&= 
\big|\delta(A_1\cap A_2)\big|+\big|\delta(A_2\setminus A_1)\big|+ \big|\delta(A_1\setminus A_2)\big|+\big|\delta({A_1\cup A_2})\big|.
\end{split}
\end{equation}

We begin with Frank's version of Mader's lifting theorem.

\begin{thm}[{Frank 1992 \cite{frank}}]
\label{frank} Let $G$ be a {finite} graph containing a vertex $s$ not incident with a cut-edge of $G$ such that $\deg(s)\neq 3$. Then there are $\lfloor \deg(s)/2 \rfloor $ pairwise disjoint {$\lambda_G$}-admissible pairs of edges incident with $s$. 
\end{thm}

{Recall that if a pair of edges is $\lambda_G$-admissible, then it is $\tau$-admissible for any target function $\tau$. Hence, } the same conclusion {in Theorem~\ref{frank}} holds for any target function $\tau$. A consequence of Frank's theorem is that the size of any set of edges incident $s$ such that every two edges in the set do not form a $\tau$-admissible pair is at most $\lceil \deg(s)/2 \rceil$. In particular, Frank's theorem implies the same upper bound on the maximum size of an independent set in $L(G,s,\tau_A)$.

Let $\tau$ be any target function. 
Let $D$ be a subset of $V(G)-s$. Then $r_{\tau}(D)$ is defined to be $\max\set{\tau(a,b)}:{a \in D, b \notin D \cup \{s\}}$. We say the set $D$ is $\tau$-\emph{dangerous} if $|\delta_G(D)| \leq r_{\tau}(D) + 1$.

\begin{thm}[Theorem 1.1 in \cite{ok2016liftings}]
\label{dangerous set theorem}
Let $G$ be a finite graph and let $s$ be a vertex of $G$ that is not incident with a cut-edge, and $\deg(s)\neq 3$. Let $F$ be a set of at least two edges, all incident with $s$. Then no pair of edges in $F$ defines a $\tau$-admissible lift if and only if there is a $\tau$-dangerous set $D$ so that, for every $sv\in F$, $v\in D$.
\end{thm}

{Note that a set $F$ of edges incident with $s$ such that no pair of edges in $F$ is $\tau$-admissible is an independent set of vertices in $L(G,s,\tau)$. A dangerous set $D$ as given by the theorem is called a \emph{corresponding} dangerous set to the independent set $F$.}

In the following, whenever we speak about the lifting graph $L(G,s,\tau_A)$, we assume implicitly that we have the hypotheses of Theorem~\ref{main lifting theorem}.

\textbf{Throughout the rest of the section, we have the following assumptions:}

\begin{itemize}
\item $G$ is a finite graph, and $\ell \geq 2$ is an integer, 
\item $A \subset V(G)$ is a proper subset such that $\lambda_G(x,y)\geq \ell$ for any two vertices $x$ and $y$ in $A$,
\item $s \in V(G) \setminus A$ is such that $\deg(s)>3$,
\item there is no edge with both end-vertices outside of $A$, {so in particular all the neighbours of $s$ are in $A$,} and since $\ell \geq 2$ this also implies that $s$ is not incident with a cut-edge (as assumed in Theorem \ref{dangerous set theorem} and Frank's Thoerem \ref{frank}), and finally,
\item $\tau_A$ assigns value $\ell$ to pairs in $A$, and $0$ otherwise. Later, we will set $\ell = 2k$, but for the moment $\ell$ might be odd, too.
\end{itemize}
In particular, we apply Theorem~\ref{dangerous set theorem} always with our specific target function $\tau_A$ in mind.
Then for a set $D$ of vertices being \textbf{dangerous} means,
\begin{center}
 both $D$ and $V\setminus D$ contain vertices from $A$ and $|\delta_G(D)|\leq \ell+1$. 
\end{center}

\begin{lemma}[Adaptation of Lemma 3.2 from {\cite{ok2016liftings}}]
 \label{upper bound on two intersecting dangerous sets}
Let $F_1$ and {$F_2$} be two independent sets in $L(G,s,\tau_A)$ of size $r_1$ and $r_2$ respectively, and suppose there are dangerous sets $D_1$ and $D_2$ so that {$F_i= \delta(s)\cap \delta(D_i)$ for $i =  1,2$}. Set $\alpha=|F_1\cap F_2|$. If $ \alpha > 0$, $r_1$ and $r_2$ are greater than $\alpha$ and $\overline{D_1\cup D_2}$ {meets} $A$, then $r_1+r_2\leq\lfloor{\deg(s)/2}\rfloor + 2 $.
\end{lemma}
\begin{proof}
Since $D_1$ and $D_2$ are dangerous, we have $|\delta_{G-s}(D_i)|\leq (\ell+1)-r_i$ for $i\in \{1,2\}$. By the hypotheses of the lemma, and since the neighbours of $s$ are all in $A$, then each one of $D_1\cap D_2$, $D_1\setminus D_2$, $D_2\setminus D_1$, and $D_1\cup D_2$ separate two vertices from $A$. Thus, by the connectivity condition on $A$, $|\delta_{G-s}(D_1\cap D_2)|\geq \ell-\alpha$, $|\delta_{G-s}(D_1\setminus D_2)|\geq \ell-(r_1 - \alpha)$, $|\delta_{G-s}(D_2\setminus D_2)|\geq \ell-(r_2 - \alpha)$, and $|\delta_{G-s} (D_1\cup D_2)|{= |\delta_{G-s} (\overline{D_1\cup D_2})|}\geq \ell - (\deg(s) - (r_1+r_2-\alpha))$.

Together with equation (\ref{standard equation}), these inequalities imply
$$2(\ell+1-r_1+\ell+1-r_2)\geq (\ell-\alpha)+(\ell-(r_2-\alpha))+(\ell-(r_1-\alpha))+ (\ell - (\deg(s) - (r_1+r_2-\alpha))).$$
This yields $\deg(s)+4\geq 2 (r_1+r_2)$. Thus, $\lfloor{\deg(s)/2}\rfloor + 2 \geq r_1+r_2$ as required.
\end{proof}

The following proposition proves Theorem \ref{main lifting theorem} for $\deg(s)=4$.

 \begin{prop}[Adaptation of {\cite[Proposition~3.4]{ok2016liftings}}]
 \label{PropDeg(s)=4}
If $\deg(s)=4$, then $L(G,s,\tau_A)$ is one of: $K_4$, $K_{2,2}$, or the disjoint union of two edges. The last case can hold only if {$\ell$} is odd.
 \end{prop}
\begin{proof}
The proof is almost exactly the same as in \cite{ok2016liftings}. Let $e_1, e_2, e_3, e_4$ be the edges incident with $s$. Note that if a pair is {admissible}, then so is the complementary pair. It follows {from Theorem~\ref{frank}} that $L(G,s,\tau_A)$ is the union of perfect matchings. Thus, it is one of $K_4$, $K_{2,2}$, or the disjoint union of two edges. In the first two cases it is a complete multipartite graph. We now show that the last case cannot happen if $\ell$ is even.

If $L(G,s,\tau_A)$ is the disjoint union of two edges, then one of the edges incident with $s$, say $e_1$ is not liftable with two of the other edges incident with $s$, say $e_2$ and $e_3$.
By Theorem~\ref{dangerous set theorem}, for {each} $i\in \{2,3\}$, there {is a } dangerous set $D_i$ containing the non-$s$ {endvertices} of $e_1$ and $e_i$. Since $L(G,s,\tau_A)$ is the disjoint union of two edges, it does not contain an independent set of size $3$. Thus, $e_3$ and $e_4$ have their non-$s$ {endvertices} in $\overline{D_2 \cup \{s\}}$ and, $e_2$ and $e_4$ have their non-$s$ {endvertices} in $\overline{D_3 \cup \{s\}}$. In particular $e_4$ has its non-$s$ {endvertices} in $\overline{D_2\cup D_3 \cup \{s\}}$.

  Each one of $D_2\cap D_3$, $D_2\setminus D_3$, $D_3\setminus D_2$, and $\overline{D_2\cup D_3 \cup \{s\}}$ contains exactly one neighbour of $s$, so each one of these sets separates $A$ by our standing assumption on $A$. Thus, each one of $|\delta_{G-s} (D_2\cap D_3)|$, $|\delta_{G-s} (D_2\setminus D_3)|$, $|\delta_{G-s} (D_3\setminus D_2)|$, and $|\delta_{G-s}(\overline{D_2\cup D_3 })|$ is at least $\ell-1$. Also, by the definition of being dangerous, we have $|\delta_G(D_i)|\leq \ell+1$, so $|\delta_{G-s}(D_i)| \leq \ell-1$. 
  
  For equation (\ref{standard equation}) to hold, all these inequalities have to hold with equality, and $|E_{G-s}(D_2\cap D_3, \overline{D_2\cup D_3})|=|E_{G-s}(D_2\setminus D_3, D_3\setminus D_2)|=0$.
  But then it is a routine calculation to verify that 
  $$|E_{G-s}(D_3\cap D_2, D_3\setminus D_2)| = |E_{G-s}(D_3\cap D_2, D_2\setminus D_3)|,$$
   from which it follows that $|\delta_{G-s}(D_3\cap D_2)| = \ell-1$ must be even, so $\ell$ is odd.
  \end{proof}

 \begin{prop}[Adaptation of {\cite[Proposition~3.1]{assem2022analysis}}]
 \label{prop_odddegree} 
Assume that $L(G,s,\tau_A)$ contains an independent set of size $\lceil{\deg(s)/2}\rceil$. If $\deg(s)$ is odd and at least $5$, then $L(G,s,\tau_A)$ is either
                  \begin{enumerate}
                  \item[(i)] {a complete bipartite graph with partition classes differing in size by exactly $1$, or} 
                  
                  \item[(ii)] {the union of an isolated vertex and a complete bipartite graph with equally sized partition classes.} 
                      
                  \end{enumerate}
\end{prop}
\begin{proof}
This is almost exactly the same as the proof of Case 1 in \cite[Proposition~3.5]{ok2016liftings}. Let $F$ be an independent in $L(G,s,\tau_A)$ of size $\lceil \deg(s)/2 \rceil$. By Theorem~\ref{dangerous set theorem} there is a dangerous set $D_1$ containing the non-$s$ {endvertices} of the edges in $F$. Then $|\delta_G(D_1)|\leq \ell+1$.

 By Frank's Theorem \ref{frank}, $F$ has the maximum possible size of an independent set in $L(G,s,\tau_A)$, therefore, 
 {$E(s,D_1)=F$}.
 Thus, $|{E(s,D_1)}|=|F|=\p{\deg(s)+1}/2$, and $|E(s,\overline{D_1\cup \{s\}})| = \p{\deg(s)-1}/2$.
 Note that the disjoint unions of $E(D_1,\overline{D_1\cup \{s\}})$ with $E(s,D_1)$ and with $E(s,\overline{D_1\cup \{s\}})$ respectively give $\delta(D_1)$ and $\delta(\overline{D_1\cup \{s\}})$. Thus, $|\delta(\overline{D_1\cup \{s\}})|=|\delta(D_1)|-1\leq \ell$. 
 It follows that $\overline{D_1\cup \{s\}}$ is dangerous, too, as it has small boundary and separates vertices from $A$ as $D_1$ does. Therefore, also the complement $\overline{F}$ is an independent set in $L(G,s,\tau_A)$. Now we show that 
\begin{enumerate}
\item[(i)] either every edge in $F$ lifts with every edge in $\overline{F}$ (this gives the complete bipartite possibility) or 
\item[(ii)] there is a unique edge in $F$ that does not lift with any edge in $\delta(s)$, but any other edge of $F$ lifts with each edge of $\overline{F}$. This gives the isolated vertex plus complete bipartite case.
\end{enumerate}

More precisely, we show that if an edge $e_1\in F$ does not lift with an edge $e_2\in\overline{F}$, then $e_1$ does not lift with any other edge, and then by Frank's Theorem \ref{frank} it is the only such edge. We present this in the following claim.

\begin{clm}\label{claim odd}
Let $e_1\in F$ and $e_2\in \overline{F}$ be a pair of non-admissible edges. Then $e_1$ is isolated in $L(G,s,\tau_A)$.
\end{clm}
\begin{proof}
Let $D_2$ be a dangerous set that contains the non-$s$ endvertices of $e_1$ and $e_2$. The maximum possible size of an independent set in $L(G,s,\tau_A)$ is $\lceil \deg(s)/2 \rceil$; consequently, at most $\lceil \deg(s)/2 \rceil$ of the edges incident with $s$ have their non-$s$ {endvertices} in $D_2$. Therefore, at least $\lfloor{\deg(s)/2}\rfloor$ of the edges incident with $s$ have their non-$s$ ends in $\overline{D_2\cup \{s\}}$. Now since $e_2$ is in $\overline{F}$ but has its non-$s$ end in $D_2$ and $|\overline{F}|=\lfloor{\deg(s)/2}\rfloor$, the set of edges incident with $s$ whose non-$s$ {endvertices} are in $\overline{D_2\cup \{s\}}$ contains an edge $e$ from $F$. Also since $e_1\in F$ has its non-$s$ {endvertex} in $D_2$, $e$ is in $F\setminus \{e_1\}$.

For $i=1,2$, set $F_i=\delta(D_i)\cap \delta(s)$; so $F_1=F$. Three of the hypotheses of Lemma~\ref{upper bound on two intersecting dangerous sets} are satisfied: $e_1\in F_1\cap F_2$ ($\alpha>0$), $e_2\in F_2\setminus F_1$ ($|F_2|>\alpha$), and $e\in F_1\setminus F_2$ ($|F_1|>\alpha$). If the other hypothesis of the lemma, that $\overline{D_1\cup D_2 \cup \{s\}}$ contains vertices from $A$, is also satisfied, then $|F_1|+|F_2|\leq \lfloor \deg(s)/2 \rfloor +2$, i.e.
$\lceil \deg(s)/2 \rceil+ |F_2|\leq \lfloor{\deg(s)/2}\rfloor + 2 $. Since $\deg(s)$ is odd, this means that $|F_2|\leq 1$, a contradiction since $F_2$ contains both $e_1$ and $e_2$. Thus, $\overline{D_1\cup D_2 \cup \{s\}}$ does not contain vertices from $A$, so all the edges incident with $s$ have their other {endvertices} in $D_1\cup D_2$ (as the neighbours of $s$ belong to $A$ by our standing assumption). In particular, all edges from $\overline{F}$ have their non-$s$ {endvertices} in $D_2\setminus D_1$. 

Now since $|\delta(s)\setminus F|= \lfloor{\deg(s)/2}\rfloor$ and since there is an edge from $s$ to $D_1\cap D_2$, we have $|\delta(D_2)\cap \delta(s)|\geq \lceil \deg(s)/2 \rceil$. This is the maximum possible size of an independent set in $L(G,s,\tau_A)$, and since $D_2$ is dangerous, $|\delta(D_2)\cap \delta(s)|= \lceil \deg(s)/2 \rceil$. 

Consequently, $|E(s, D_1 \cap D_2)|=1$, and all the non-$s$ {endvertices} of the edges of $F$ other than $e_1$ are in $D_1\setminus D_2$. Since all the non-$s$ {endvertices} of the edges incident with $s$ other than $e_1$ are either in $D_1\setminus D_2$ or $D_2\setminus D_1$, and $e_1$ has its non-$s$ {endvertex} in $D_1\cap D_2$, $e_1$ does not lift with any other edge. 
\end{proof}

By Frank's Theorem \ref{frank}, there can only be one edge $e$ in $\delta(s)$ that is not liftable with any other edge in $\delta(s)$. If such an edge $e$ exists, and if $f$ is an edge in $F\setminus\{e\}$ and $f'$ is an edge in $\delta(s)\setminus F$, then $f$ and $f'$ is an $\tau_A$-admissible pair, otherwise, by Claim \ref{claim odd}, $f$ is not liftable with any edge in $\delta(s)\setminus \{f\}$, contradicting the uniqueness of $e$. Thus every edge in $F\setminus \{e\}$ lifts with every edge in $\overline{F}$. This gives the structure of an isolated vertex plus a balanced complete bipartite graph for the lifting graph. 
\end{proof}

\begin{prop}[Adaptation of {\cite[Lemma~3.3]{assem2022analysis}}]
\label{prop_maximalindependentsets}
 
 Assume that $\deg(s) \geq 4$ and $I_1$ and $I_2$ are two distinct maximal independent sets in $L(G,s,\tau_A)$ each of size at least $2$.

\begin{enumerate}

\item[(1)] Then $|I_1\cap I_2|\leq 1$.

\item [(2)] If $I_1\cap I_2\neq \emptyset$ and $I_1\cup I_2 \neq V(L(G,s,\tau_A))$, then $\ell$ is odd.

\end{enumerate}

\end{prop}

\begin{proof}
The proof follows a generalized version of the proof of Case 2 in \cite[Proposition~3.5]{ok2016liftings}. {First note that since $\deg(s)>3$, the maximum possible size of an independent set is $\lceil \deg(s)/2 \rceil$ by Frank's Theorem \ref{frank}. Therefore, if $I_1\cap I_2 \neq \emptyset$ and $I_1\cup I_2$ is the entire vertex set of $L(G,s,\tau_A)$, then $I_1$ and $I_2$ are both of size $\lceil \deg(s)/2 \rceil$ and they intersect in exactly one vertex, otherwise they will not cover $V(L(G,s,\tau_A))$. So this yields assertion (1) of the proposition in the case where $I_1\cup I_2 = V(L(G,s,\tau_A))$. Thus, to prove (1) and (2), we may assume for the remainder of the proof that $I_1\cap I_2\neq \emptyset$ and that $I_1\cup I_2 \neq V(L(G,s,\tau_A))$.}

{Let $D_1$ and $D_2$ be two dangerous sets in $G$ corresponding to $I_1$ and $I_2$ respectively.} Let $k_1= |E(\{s\}, D_1\setminus D_2)|$, $k_2=|E(\{s\}, D_2\setminus D_1)|$, $k_3=|E(\{s\}, D_1\cap D_2)|$, then $k_3\neq 0$.

Since $D_1$ and $D_2$ are dangerous, we have $|\delta_G(D_1)|\leq \ell+1$ and $|\delta_G(D_2)|\leq \ell+1$. In particular,
$|\delta_{G-s}(D_1)|\leq (\ell+1)-(k_1+k_3)$, and
$|\delta_{G-s}(D_2)|\leq (\ell+1)-(k_2+k_3)$.

Since $I_1$ and $I_2$ are distinct maximal independent sets in $L(G,s,\tau_A)$, {the union} $I_1\cup I_2$ is not independent, consequently, $D_1\cup D_2$ is not dangerous.
Also because $I_i$ is a maximal independent set of $L(G,s,\tau_A)$ for $i\in \{1,2\}$, each $D_i$ does not contain neighbours of $s$ other than those that are end-vertices of edges in $I_i$.

This and the assumption that $I_1\cup I_2 \neq V(L(G,s,\tau_A))$ ({because $I_1\cap I_2 \neq \emptyset$}) imply that at least one edge incident with $s$ has its non-$s$ end in the set $\overline{D_1\cup D_2\cup \{s\}}$. Thus, $D_1\cup D_2$ separates vertices from $A$ in $G-s$. 
By the definition of a dangerous set, the only other way for $D_1\cup D_2$ to not be dangerous is if $|\delta_G(D_1\cup D_2)|\geq \ell+2$. 
  This means that $|\delta_{G-s}(D_1\cup D_2)|\geq (\ell+2)-(k_1+k_2+k_3)$.

Since $D_1\cap D_2$, $D_1\setminus D_2$, and $D_2\setminus D_1$ also separate $A$, we have, $|\delta_{G-s}(D_1\cap D_2)|\geq \ell-k_3$, $|\delta_{G-s}(D_1 \setminus D_2)|\geq \ell-k_1$, and $|\delta_{G-s}(D_2 \setminus D_1)|\geq \ell-k_2$. Observe that in $G-s$:

\begin{align*}
\begin{split}
& 2\bigg[\big|\delta(D_1)\big| +\big|\delta(D_2)\big|- \Big(\big|E(D_1 \cap D_2, \overline{D_1 \cup D_2})|+\big|E(D_2\setminus D_1, D_1\setminus D_2)\big|\Big)\bigg]
\\
&\leq 2\bigg [(\ell+1)-(k_1+k_3)+ (\ell+1)-(k_2+k_3) \bigg]
\\
& =4\ell-2(k_1+k_2+k_3)+2+(2-2k_3)
\end{split}
\end{align*}
and  
\begin{align*}
\begin{split}
&\big|\delta (D_1\cap D_2)\big|+\big|\delta (D_2\setminus D_1)\big|+ \big |\delta (D_1\setminus D_2)\big|+ \big|\delta (D_1\cup D_2)\big| \geq \\ & (\ell-k_3)+(\ell-k_2)+(\ell-k_1)+(\ell+2)-(k_1+k_2+k_3)\\
& =4\ell-2(k_1+k_2+k_3)+2.
\end{split}
\end{align*}
By Equation~(\ref{standard equation}), it follows that $2-2k_3\geq 0$, i.e. $k_3\leq 1$ as desired, {yielding assertion (1) of the proposition.}

If $k_3=1$, then the inequalities above have to be equalities throughout. More precisely,
\begin{enumerate}
\item[{(i)}] $|\delta_{G-s}(D_1)|=|\delta_{G-s}(D_1 \setminus   D_2)|= \ell-k_1$, so $|\delta_{G}(D_1\setminus D_2)|=\ell$ and $|\delta_G(D_1)|=\ell+1$,
\item[{(ii)}] $|\delta_{G-s}(D_2)|= |\delta_{G-s}(D_2 \setminus D_1)|=\ell-k_2$, so $|\delta_{G}(D_2\setminus D_1)|=\ell$ and $|\delta_G(D_2)|=\ell+1$,
\item[{(iii)}] $|\delta_{G-s}(D_1\cap D_2)|= \ell-1$,
\item[{(iv)}] $|\delta_{G-s}(D_1\cup D_2)|=(\ell+2)-(k_1+k_2+k_3)=\ell+1-k_1-k_2$, and
\item[{(v)}] $|E_{G-s}(D_1\cap D_2, \overline{D_1\cup D_2})|=|E_{G-s}(D_2\setminus D_1, D_1\setminus D_2)|=0$.
\end{enumerate}
 {From (v) it follows that} $|\delta_{G-s}(D_2)|=|E_{G-s}(D_2\setminus D_1, \overline{D_1\cup D_2})|+|E_{G-s}(D_1\cap D_2, D_1\setminus D_2)|$, and \newline
$|\delta_{G-s}(D_2\setminus D_1)|=|E_{G-s}(D_2\setminus D_1,\overline{D_1\cup D_2})| +|E_{G-s}(D_2\setminus D_1, D_1\cap D_2)|$. {From (ii) we have} $|\delta_{G-s}(D_2)|= |\delta_{G-s}(D_2 \setminus D_1)|$. Cancelling the common $|E_{G-s}(D_2\setminus D_1, \overline{D_1\cup D_2})|$ on both sides yields $$|E_{G-s}(D_2\setminus D_1, D_1\cap D_2)|=|E_{G-s}(D_1\cap D_2, D_1\setminus D_2)|.$$
Now this last equality, {(iii), and (v)} imply that 
\begin{align*}
\ell-1 & =|\delta_{G-s}(D_1\cap D_2)| \\ 
& = |E_{G-s}(D_1\cap D_2, D_1\setminus D_2)|+ |E_{G-s}(D_1\cap D_2, D_2\setminus D_1)| \\
 & = 2 \cdot |E_{G-s}(D_1\cap D_2, D_2\setminus D_1)|.
 \end{align*}
Thus $\ell$ is odd, {yielding assertion (2) of the proposition.}
\end{proof}

We now have all ingredients for the proof of our main theorem, which we restate here for convenience:

\begin{thm} 
Let $G$ be a finite graph, {$k$ a positive} integer, and $A$ a proper subset of $V(G)$ such that between any two vertices in $A$ there are $2k$ edge-disjoint paths in $G$. Assume that there is no edge with both endvertices outside of $A$, and let $s$ be a vertex not in $A$ with $\deg(s) > 3$. Then $L(G,s,\tau_A)$ is either
\begin{itemize}

\item a complete multipartite graph, or

\item {the union of an isolated vertex and a complete bipartite graph with equally sized partition classes. Hence, this latter case can only occur if $\deg(s)$ is odd.}

\end{itemize}

\end{thm}

\begin{proof}
{If $\deg(s) = 4$, then the result follows from Proposition~\ref{PropDeg(s)=4} since $2k$ is even. So assume that $\deg(s) \geq 5$. If all maximal independent sets in $L(G,s,\tau_A)$ are pairwise disjoint, then  $L(G,s,\tau_A)$ is a complete multipartite graph. Otherwise, there are two distinct intersecting maximal independent sets $I_1$ and $I_2$. By Proposition~\ref{prop_maximalindependentsets}, since $2k$ is even, we have $I_1\cup I_2 = V(L(G,s,\tau_A))$. By Frank's Theorem \ref{frank}, $I_1$ and $I_2$ have size at most $\lceil \deg(s)/2 \rceil$. Thus, for their non-disjoint union to cover $V(L(G,s,\tau_A))$, they must both be of size $\lceil \deg(s)/2 \rceil$ with $\lceil \deg(s)/2 \rceil > \deg(s)/2$. Therefore, $\deg(s)$ is odd, and the assertion follows from Lemma~\ref{prop_odddegree}.}
\end{proof}

\section{Extending orientations of Eulerian subgraphs}

For a graph $G=(V,E)$, a set of edges $F \subset E$ and a non-trivial bipartition $(A,B)$ of $V$, we write  $F(A,B)$ for the set of edges in $F$ with one endvertex in  $A$ and the other in $B$. {If $G$ is directed, then $\vec{F}(A,B)$ denotes the set of edges in $F$ from $A$ into  $B$. A directed graph $G=(V,E)$ is \emph{balanced} if $|\vec{E}(A,B)| = |\vec{E}(B,A)|$ for all non-trivial bipartition $(A,B)$ of $V$.}

Suppose we are given a $2k$-edge-connected graph $G$. In order to define a $k$-arc-connected orientation, we need to orient the edges of any cut $E(A,B)$ in $G$ such that at least $k$ edges go from $A$ to $B$, and at least $k$ edges go from $B$ to $A$. Then clearly, if we orient some Eulerian subgraph $H$ of $G$ \emph{consistently} (that is, along some closed Euler trail of $H$), then from any cut $E(A,B)$ in $G$, there will be an equal number of oriented edges (we say the cut is \emph{balanced}) in both directions (and possibly other edges which are not yet oriented). Thus, on our way towards a $k$-arc-connected orientation for $G$, we have made no obvious mistake yet.

In fact, it is known that in any $2k$-edge-connected graph $G$, any consistent orientation of a closed, Eulerian subgraph $H$ extends to a $k$-arc-connected orientation of $G$, \cite[Theorem~9.2.3]{frank2011connections}.

As in the previous section, for two vertices $x,y$ in a graph $G$, we write $\lambda(x,y)$ for maximum number of edge-disjoint $x-y-$paths, and $\lambda^*(x,y)$ for the greatest even number not exceeding $\lambda(x,y)$. Further, for two vertices $x,y$ in an oriented graph $\vec{G}$ define $\alpha(x,y)$ as the maximum number of arc-disjoint directed paths from $x$ to $y$.
Let us say an orientation $\vec{G}$ of a graph $G$ is \emph{well-balanced}, if
$$\alpha(x,y) \geq \frac{\lambda^*(x,y)}{2}$$ for any two distinct vertices $x,y \in G$. 

Recall that the \emph{Strong orientation theorem} of Nash-Williams says that every finite graph admits a well-balanced orientation of the whole graph \cite{nash1960orientations}. Again, it is known that every consistent orientation of a closed Eulerian subgraph extends to a well-balanced orientation \cite[Corollary~2]{kiraly2006simultaneous}. As our last ingredient, we observe that Nash-Williams's original proof for the orientation theorem also yields the same extension property for all open Eulerian partial orientations.

\begin{thm}\label{thm_eulerpreorientation}
Let $G$ be a finite graph and $H \subseteq G$ an open or closed Eulerian subgraph. Then any consistent orientation $\vec H$ of $H$ can be extended to a well-balanced orientation of $G$.
\end{thm}
\begin{proof}
An \emph{odd vertex pairing} of a finite graph $G=(V,E)$ is a partition $P$ of the vertices of odd-degree in $G$ into sets of size two. Interpreting $P$ as edges, we obtain an Eulerian graph $G' = (V,E')$ where $E'  = E \,  \dot\cup \, P$. Then $H \subseteq G \subseteq G'$.
Nash-Williams showed in \cite{nash1960orientations}*{Theorem~2} that every graph $G = (V,E)$ has an odd-vertex pairing $P$ such that for every two $x,y \in V$ and every bipartition $(X,Y)$ of $V$ with $x \in X$ and $y \in Y$ holds: 
$$(\star) \quad \quad |E(X,Y)| - |P(X,Y)| \geq \lambda^*(x,y).$$

We claim that with such an odd-vertex pairing, any consistent orientation $\vec G'$ of the Eulerian graph $G'$ that extends $\vec{H}$ restricts to a well-balanced orientation $\vec{G}$ of $G$ as desired.  

To see this claim, fix any two vertices $x$ and $y$, and let $(X,Y)$ be any partition of $V(G)$ separating $x$ from $y$. Since $\vec G'$ is balanced, it follows that
$$|\vec{E}(X,Y)| + |\vec{P}(X,Y)|  =\frac{|E(X,Y)| + |P(X,Y)|}{2}.$$
However, since $|\vec{P}(X,Y)|  \leq  |P(X,Y)|$, it follows that
\begin{equation*}
 |\vec{E}(X,Y)|  \geq \frac{|E(X,Y)| + |P(X,Y)|}{2} - |P(X,Y)| = \frac{|E(X,Y)| - |P(X,Y)|}{2}  \, \stackrel{(\star)}{\geq} \, \frac{\lambda^*(x,y)}{2} 
\end{equation*}
and hence $ \alpha(x,y) \geq \lambda^*(x,y)/2 $ as desired.
\end{proof}

\section{Nash-Williams' orientation theorem for infinite graphs}
\label{sec_orientation}

In this section, we are finally ready to prove our orientation results. Let us first remark that by substituting Theorem~\ref{thm_eulerpreorientation} for \cite[Theorem~6]{thomassen2016orientations} in Thomassen's proof improves Thomassen's $8k$ bound to $4k$ for all graphs, {improving on \cite{assem2023towards} where this conclusion has been shown for one-ended, locally finite graphs. Hence, in the remainder of this section we focus }on the optimal $2k$ bound for locally finite graphs with at most countably many ends.

\begin{thm}
Let $G$ be a locally finite graph with countably many ends, and $k$ a positive integer. Then $G$ is $2k$-edge-connected if and only if it has a $k$-arc-connected orientation.
\end{thm}

\begin{proof} 
By Egyed's result, we may assume that $k \geq 2$. It is clear that $2k$-edge-connectivity is a necessary condition for the existence of a $k$-arc-connected orientation, so here we prove the other implication. Enumerate $V=\{v_0,v_1,...\}$. Beginning with $A_0 = \{ v_0 \}$ and {$\vec W_0$ any consistent orientation of a cycle $W_0 \subseteq G$} containing $v_0$, we will construct a sequence of finite, $2$-edge-connected subgraphs $W_0 \subseteq W_1 \subseteq W_2 \subseteq \cdots$ of $G$ with compatible orientations $\vec W_0 \subseteq \vec W_1 \subseteq \vec W_2 \subseteq \cdots$ and sets of vertices $A_0 \subseteq A_1 \subseteq A_2 \cdots$ such that for all $n\geq 0$: 
\begin{enumerate}[label=(\roman*)]
    \item\label{item_exhaust} $\{v_0,\ldots,v_n\} \subseteq A_n \subseteq V(W_n)$.
    \item\label{item_Eulerian} For every component $B$ of $G-A_n$, {of all vertices in $V(B) \cap V(W_n)$} all but possibly at most one exceptional vertex  have in-degree equalling out-degree in $\vec W_n$, with the exceptional vertex having a difference of $1$ between in- and out-degree.
    \item\label{item_kconnect} For every two distinct vertices $x,y$ in $A_n$, there are $k$ arc-disjoint directed paths in $\vec W_n$ from $x$ to $y$ and from $y$ to $x$.
\end{enumerate}
Once the construction is complete, we claim that properties \ref{item_exhaust} and \ref{item_kconnect} imply that any orientation $\vec G$ of $G$ extending $\vec{W}:=\bigcup_{i \in \mathbb{N}} \vec W_i$ is $k$-arc-connected. Indeed, for every two distinct vertices $x,y$ in $G$, by \ref{item_exhaust} there is an $i \in \mathbb{N}$ with $x,y \in A_i$, and so by \ref{item_kconnect} there are $k$ arc-disjoint directed paths in $\vec W_i$ from $x$ to $y$ and from $y$ to $x$. Since $\vec W_i \subseteq \vec W$ as oriented subgraphs, these directed paths are directed also in $\vec W$, and hence in $\vec G$, as desired.

Thus, it remains to describe the inductive construction, and this is where property \ref{item_Eulerian} is needed. So suppose inductively that we have already constructed $A_n$ and $\vec W_n$ according to \ref{item_exhaust}--\ref{item_kconnect}. Since $G$ has countably many ends, we may apply Theorem~\ref{finite partition} to the set $A'_{n+1} := V(W_n) \cup \{v_{n+1}\}$ to obtain a finite set $A_{n+1} \supseteq A'_{n+1}$ such that the components of $G-A_{n+1}$ are boundary-linked sets. Applying Theorem~\ref{immersion} yields a finite set $X$ meeting every boundary-linked component at most once, and also an immersion $W_{n+1}$ on $A_{n+1} \cup X$ in $G$ of a finite, $3$-edge-connected graph $H$ (with $d_{H}(x)=3$ for all $x \in X$ and $\lambda_{H} (a,b) \geq 2k$ for all  $a,b \in A_{n+1}$).
By (iii) in Theorem~\ref{immersion}, we have $W_n \subseteq H$.
\begin{figure}
    \centering
\begin{tikzpicture}
\tikzstyle{O}=[draw,very thick,black!30]
        \tikzstyle{dot}=[draw,circle,color=white,minimum size=0pt,inner sep=0pt]
                            \tikzstyle{edge} = [draw,line width=0.5pt,-]
    \tikzstyle{green node}=[draw,circle,fill=black,minimum size=4pt,
                            inner sep=0pt]
  \tikzstyle{edge} = [draw,line width=0.6pt,-]

    \tikzstyle{black node}=[draw,circle,fill=black,minimum size=3pt,
                            inner sep=0pt]

\draw[thick] (-4.7,3.8) .. controls (-2.25,-0.3) and (-0.75,-0.3) .. (1.2,3.8);

\draw[thick] (-4.5,3.8) .. controls (-3.6,2.3) and (-2.1,2.3) .. (-2.7,3.8);
\draw[thick] (-2.5,3.8) .. controls (-2.1,2.3) and (-0.6,2.3) .. (-0.9,3.8);
\draw[thick] (-0.7,3.8) .. controls (-0.6,2.3) and (0.7,2.3) .. (0.9,3.8);

\draw[thick,color=blue] (0.2,2.1) .. controls (-1.1,2.7) and (-2.2,2.7) .. (-3.5,2.1);

\draw[thick,color=blue] (-4.7,0.9) .. controls (-5,0.5) and (-5.2,-0.5) .. (-5,-2.3);
\draw[thick] (-6,2) .. controls (-2.5,-0.6) and (-2.5,-1.2) .. (-6,-3);
\draw[thick] (-6,1.8) .. controls (-4.7,1) and (-4.7,0.5) .. (-6,0);

\draw[thick] (2,2.7) .. controls (1,0.5) and (1.2,-1) .. (3,-1.3);

\draw[thick] (-6,-0.8) .. controls (-4.9,-1) and (-4.9,-2.5) .. (-6,-2.6);
                           \draw (-3.3,2) node[dot] (j1) []{};
                            \draw (-3.1,2.1) node[dot] (j2) []{};
                            \draw (-2.9,2.2) node[dot] (j3) []{};
                            \draw (-2.7,2.25) node[dot] (j4) []{};

\draw[thick] (-6,1.8) .. controls (-4.7,1) and (-4.7,0.5) .. (-6,0);
                           \draw (-4.65,0.65) node[dot] (k1) []{};
                            \draw (-4.7,0.5) node[dot] (k2) []{};
                            \draw (-4.75,0.35) node[dot] (k3) []{};
                            \draw (-4.8,0.2) node[dot] (k4) []{};

\draw[rounded corners, thick] (k1)  -- ++(160:1.2cm) -- ++(220:0.3cm) -- ++(270:0.35cm) -- ++(k3);

\draw[rounded corners, thick] (k2)  -- ++(160:1.4cm) -- ++(220:0.3cm) -- ++(270:0.55cm) -- ++(k4);

                           \draw (-4.9,-2) node[dot] (l1) []{};
                            \draw (-4.9,-1.8) node[dot] (l2) []{};
                            \draw (-4.9,-1.6) node[dot] (l3) []{};
                            \draw (-4.95,-1.4) node[dot] (l4) []{};
\draw[rounded corners, thick] (l2)  -- ++(180:0.7cm) -- ++(220:0.3cm) -- ++(270:0.35cm) -- ++(l1);

\draw[rounded corners, thick] (l4)  -- ++(160:1.2cm) -- ++(220:0.3cm) -- ++(270:0.55cm) -- ++(l3);
                            
\draw[thick,color=blue] (1.7,1.8) .. controls (2.3,0.6) and (3,-0.5) .. (2.7,-1.2);

\draw[thick] (2.2,2.7) .. controls (2.2,2.4) and (1,2.1) .. (3,1.2);
\draw[thick] (3.3,0.7) .. controls (2.5,0.5) and (2.3,0.3) .. (3.3,-1);
                           \draw (1.7,1.5) node[dot] (n1) []{};
                            \draw (1.8,1.2) node[dot] (n2) []{};
                            \draw (2.4,0) node[dot] (n3) []{};
                            \draw (2.5,-0.3) node[dot] (n4) []{};
                            \draw (2.6,-0.6) node[dot] (n5) []{};
                            \draw (2.7,-0.8) node[dot] (n6) []{};

\draw[rounded corners, thick] (n1)  -- ++(30:1.2cm) -- ++(-70:0.3cm) -- ++(-80:0.3cm) -- ++(n2);

\draw[rounded corners, thick] (n3)  -- ++(30:0.7cm) -- ++(-70:0.3cm) -- ++(-80:0.3cm) -- ++(n5);

\draw[rounded corners, thick] (n4)  -- ++(30:1.2cm) -- ++(-80:0.2cm) -- ++(-80:0.3cm) -- ++(n6);

                           \draw (-3.3,2) node[dot] (a1) []{};
                            \draw (-3.1,2.1) node[dot] (a2) []{};
                            \draw (-2.9,2.2) node[dot] (a3) []{};
                            \draw (-2.7,2.25) node[dot] (a4) []{};
                            \draw (-2.5,2.3) node[dot] (a5) []{};
                            \draw (-3.5,3.3) node[black node] (a) []{};

                            \path[edge] (a1) to [out=120,in=250] (a);
                            \path[edge] (a2) to [out=110,in=260] (a);
                            \path[edge] (a3) to [out=110,in=280] (a);
                            \path[edge] (a4) to [out=120,in=290] (a);
                            \path[edge] (a5) to [out=120,in=80] (a);

                             \draw (-1.9,2.3) node[dot] (b1) []{};
                             \draw (-1.7,2.3) node[dot] (b2) []{};
                             \draw (-1.5,2.3) node[dot] (b3) []{};
                             \draw (-1.3,2.3) node[dot] (b4) []{};
                             
\draw[rounded corners, thick] (b1)  -- ++(100:1.2cm) -- ++(20:0.3cm) -- ++(-30:0.35cm) -- ++(-90:0.6cm) -- ++(b3);
\draw[rounded corners, thick] (b2)  -- ++(100:1.8cm) -- ++(20:0.3cm) -- ++(-30:0.35cm) -- ++(-90:0.6cm) -- ++(b4);

\draw (-0.8,2.25) node[dot] (c1) []{};
\draw (-0.6,2.2) node[dot] (c2) []{};
\draw (-0.4,2.15) node[dot] (c3) []{};
\draw (-0.2,2.1) node[dot] (c4) []{};
\draw (0,2.05) node[dot] (c5) []{};

\draw (0.5,3.8) node[black node] (c) []{};
\path[edge] (c2) to [out=100,in=260] (c);
\path[edge] (c3) to [out=100,in=280] (c);
\path[edge] (c4) to [out=110,in=290] (c);                            

\draw[rounded corners, thick] (c1)  -- ++(70:1.4cm) -- ++(20:0.2cm) -- ++(-30:0.35cm) -- ++(-90:0.6cm) -- ++(c5);

\draw[rounded corners, thick] (b1)  -- ++(100:1.2cm) -- ++(20:0.3cm) -- ++(-30:0.35cm) -- ++(-90:0.6cm) -- ++(b3);
\draw[rounded corners, thick] (b2)  -- ++(100:1.8cm) -- ++(20:0.3cm) -- ++(-30:0.35cm) -- ++(-90:0.6cm) -- ++(b4);

                            \draw (-2.5,-1) node[dot] (x1) []{};
                            \draw (-2.5,-1.3) node[dot] (x2) []{};
                            \draw (-2.5,-0.7) node[dot] (x3) []{};
                            \draw (-2.5,-0.4) node[dot] (x4) []{};

    \begin{scope}[very thick,nodes={sloped,allow upside down}]                       
   \draw[rounded corners, thick] (x2)  --pic[pos=.5]{arrow=stealth} ++(190:1.8cm) -- ++(200:0.5cm) --pic[pos=.7]{arrow=stealth} ++(100:0.95cm)-- ++(60:0.75cm) --pic[pos=.7]{arrow=stealth} ++(x4);
   \draw[rounded corners, thick] (x1)  -- pic[pos=.7]{arrow=stealth} ++(190:1.1cm) -- ++(200:0.5cm) -- ++(60:0.35cm) --pic[pos=.5]{arrow=stealth}  ++(180:0.85cm) -- ++(60:0.35cm) --pic[pos=.7]{arrow=stealth} ++(x3);
   \end{scope}
                 \draw (-2.2,1.8) node[black node] (y) []{};
                                    
                \draw (-0.6,0.1) node[dot] (y1) []{};
                \draw (-0.8,0.1) node[dot] (y2) []{};
                \draw (-1,0.1) node[dot] (y3) []{};
                \draw (-1.2,0.1) node[dot] (y4) []{};
                \draw (-1.4,0.1) node[dot] (y5) []{};
                \draw (-1.6,0.1) node[dot] (y6) []{};
                \draw (-1.8,0.1) node[dot] (y7) []{};
                                                        
        \path[edge,postaction={on each segment={mid arrow}}] (y) to (y5);
        \path[edge,postaction={on each segment={mid arrow}}] (y) to [out=250,in=100] (y6);
        \path[edge,postaction={on each segment={mid arrow}}] (y7) to [in=220,out=150] (y);       

 \begin{scope}[very thick,nodes={sloped,allow upside down}]
    \draw[rounded corners, thick] (y1)  -- pic[pos=.7]{arrow=stealth}  ++(80:1cm) -- ++(70:0.5cm) -- ++(130:0.65cm) -- pic[pos=.3]{arrow=stealth} ++(y3);
   \draw[rounded corners, thick] (y2)  -- pic[pos=.7]{arrow=stealth} ++(80:1.1cm) -- ++(200:0.5cm) -- ++(60:0.35cm) --pic[pos=.5]{arrow=stealth} ++(180:0.85cm)  -- pic[pos=.65]{arrow=stealth} ++(y4);

                \draw[rounded corners, thick] (1.8,0.3) node[black node] (z) []{}
                     -- ++(180:0.4cm)  --pic[pos=.5]{arrow=stealth reversed} ++(30:-1.2cm) node[dot] () {};
                \draw[rounded corners, thick] (z)
                    -- ++(20:-0.6cm) --pic[pos=.5]{arrow=stealth} ++(35:-1cm) node[dot] () {};
                \draw[rounded corners, thick] (z)
                     -- ++(22:0.5cm) -- ++(-82:0.3cm) --pic[pos=.5]{arrow=stealth} ++(30:-1.7cm) node[dot] () {};
\end{scope}
 \draw[rotate around={-25:(0.9,2)},very thick,black!80] (0.8,-1.7) ellipse (60pt and 51pt);
 
 \draw (0,-1.5) node[dot] () [label=-90:$A_n$]{};

  \draw (3,-1.5) node[dot] () [label=-90:\textcolor{blue}{$A_{n+1}$}]{};
 
 \end{tikzpicture}
    \caption{An immersion on $A_{n+1} \cup X$ encompassing all previously oriented edges.}
    \label{fig:extending1}
\end{figure}
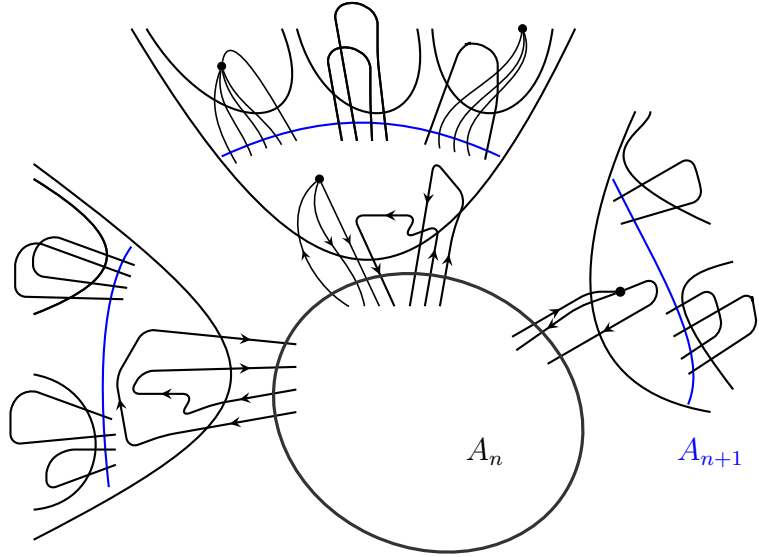

Now contract $A_n$ in $H$ to a dummy vertex $v$, and call the resulting graph $\tilde{H}$. This graph inherits from $H$ the property of being $2k$-edge-connected between the vertices of $V(\tilde{H}) \setminus X$ (including $v$) since the contracted set $A_n$ is disjoint from $X$. For each component $B$ of $G - A_n$, let $\tilde{H} \restriction B$ be the subgraph of $\tilde{H}$ induced by the dummy vertex $v$ together with $V(B)\cap V(H)$. Since $v$ is a cut-vertex in $\tilde{H}$, each $\tilde{H} \restriction B$ has the same edge-connectivity as $\tilde{H}$ between the vertices not in $X$.

Property \ref{item_Eulerian} implies that all the arcs of $\vec W_n$ {with at least on endvertex in} $B$ form a consistently oriented (open or closed) Eulerian subgraph of $\tilde{H} \restriction B$. Note that if $B$ contains an exceptional vertex $x$ as described in \ref{item_Eulerian}, then there should be a difference of $1$ between the in- and out-degree of $v$ in the oriented subgraph of $\tilde{H} \restriction B$, so $x$ and $v$ are the endvertices of an open Eulerian tour.
We can now apply Theorem~\ref{thm_eulerpreorientation} to each $\tilde{H} \restriction B$ to extend the orientation of this subgraph to a well-balanced orientation of all of $\tilde{H} \restriction B$, making this graph $k$-arc-connected between vertices of $V(\tilde{H} \restriction B) \setminus X$.
After doing this for every component $B$ of $G-A_n$, we obtain an orientation $\vec H$ of $H$.

\begin{figure}
    \centering
 \begin{tikzpicture}
\tikzstyle{O}=[draw,very thick,black!30]
        \tikzstyle{dot}=[draw,circle,color=white,minimum size=0pt,inner sep=0pt]
                            \tikzstyle{edge} = [draw,line width=0.5pt,-]
    \tikzstyle{green node}=[draw,circle,fill=black,minimum size=4pt,
                            inner sep=0pt]
  \tikzstyle{edge} = [draw,line width=0.6pt,-]

    \tikzstyle{black node}=[draw,circle,fill=black,minimum size=3pt,
                            inner sep=0pt]

\draw[thick] (-4.7,3.8) .. controls (-2.25,-0.3) and (-0.75,-0.3) .. (1.2,3.8);

\draw[thick] (-4.5,3.8) .. controls (-3.6,2.3) and (-2.1,2.3) .. (-2.7,3.8);
\draw[thick] (-2.5,3.8) .. controls (-2.1,2.3) and (-0.6,2.3) .. (-0.9,3.8);
\draw[thick] (-0.7,3.8) .. controls (-0.6,2.3) and (0.7,2.3) .. (0.9,3.8);

\draw[thick,color=blue] (0.2,2.1) .. controls (-1.1,2.7) and (-2.2,2.7) .. (-3.5,2.1);

                           \draw (-3.3,2) node[dot] (a1) []{};
                            \draw (-3.1,2.1) node[dot] (a2) []{};
                            \draw (-2.9,2.2) node[dot] (a3) []{};
                            \draw (-2.7,2.25) node[dot] (a4) []{};
                            \draw (-2.5,2.3) node[dot] (a5) []{};
                            \draw (-3.5,3.3) node[black node] (a) []{};

                            \path[edge] (a1) to [out=120,in=250] (a);
                            \path[edge] (a2) to [out=110,in=260] (a);
                            \path[edge] (a3) to [out=110,in=280] (a);
                            \path[edge] (a4) to [out=120,in=290] (a);
                            \path[edge] (a5) to [out=120,in=80] (a);

                             \draw (-1.9,2.3) node[dot] (b1) []{};
                             \draw (-1.7,2.3) node[dot] (b2) []{};
                             \draw (-1.5,2.3) node[dot] (b3) []{};
                             \draw (-1.3,2.3) node[dot] (b4) []{};
                             
\draw[rounded corners, thick] (b1)  -- ++(100:1.2cm) -- ++(20:0.3cm) -- ++(-30:0.35cm) -- ++(-90:0.6cm) -- ++(b3);
\draw[rounded corners, thick] (b2)  -- ++(100:1.8cm) -- ++(20:0.3cm) -- ++(-30:0.35cm) -- ++(-90:0.6cm) -- ++(b4);

\draw (-0.8,2.25) node[dot] (c1) []{};
\draw (-0.6,2.2) node[dot] (c2) []{};
\draw (-0.4,2.15) node[dot] (c3) []{};
\draw (-0.2,2.1) node[dot] (c4) []{};
\draw (0,2.05) node[dot] (c5) []{};

\draw (0.5,3.8) node[black node] (c) []{};
\path[edge] (c2) to [out=100,in=260] (c);
\path[edge] (c3) to [out=100,in=280] (c);
\path[edge] (c4) to [out=110,in=290] (c);                            

\draw[rounded corners, thick] (c1)  -- ++(70:1.4cm) -- ++(20:0.2cm) -- ++(-30:0.35cm) -- ++(-90:0.6cm) -- ++(c5);

                 \draw (-2.2,1.8) node[black node] (y) []{};
                                    
                \draw (-0.6,0.1) node[dot] (y1) []{};
                \draw (-0.8,0.1) node[dot] (y2) []{};
                \draw (-1,0.1) node[dot] (y3) []{};
                \draw (-1.2,0.1) node[dot] (y4) []{};
                \draw (-1.4,0.1) node[dot] (y5) []{};
                \draw (-1.6,0.1) node[dot] (y6) []{};
                \draw (-1.8,0.1) node[dot] (y7) []{};

                \draw (-1.1,-1) node[black node] (v) []{};
                                                        
        \path[edge, postaction={on each segment={mid arrow}}] (y) to (v);
        \path[edge,postaction={on each segment={mid arrow}}] (y) to [out=230,in=120] (v);
        \path[edge,postaction={on each segment={mid arrow}}] (v) to [in=200,out=150] (y);

 \begin{scope}[very thick,nodes={sloped,allow upside down}]
    \draw[rounded corners, thick] (v)  -- ++(40:0.4cm) -- pic{arrow=stealth}  ++(80:1.8cm) -- ++(70:0.5cm) -- pic{arrow=stealth} ++(130:0.65cm) -- pic{arrow=stealth} ++(v);
   \draw[rounded corners, thick] (v)  -- ++(60:0.5cm) -- pic{arrow=stealth} ++(80:2.1cm) -- ++(200:0.5cm) -- ++(60:0.35cm) --  pic{arrow=stealth} ++(180:0.85cm) -- pic{arrow=stealth} ++(v);
\end{scope}

 \draw (v) node[dot] () [label=-90:$v$]{};

 \end{tikzpicture}
    \caption{All previously oriented edges form an open or closed Eulerian tour in $\tilde{H} \restriction B$.}
    \label{fig:extending2}
\end{figure}
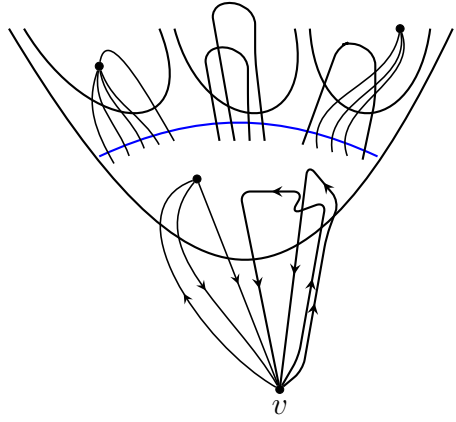

We claim that with this orientation, also $\vec H$ is $k$-arc-connected between vertices of $A_{n+1}$: Indeed, let $E(Y,Z)$ be any bond in $H$ separating two vertices of $A_{n+1}$ (that is of $V(H)\setminus X$). If $A_n$ lies completely on one side $Y$ or $Z$, then the bond restricts to a cut in some $\tilde{H} \restriction B$, and since $\tilde{H} \restriction B$ is $k$-arc-connected between the vertices of $V(\tilde{H} \restriction B) \setminus X$, there exist at least $k$ edges oriented from $Y$ to $Z$, and also from $Z$ to $Y$. And if $A_n$ meets both $Y$ and $Z$, then the bond restricts to a cut of $\vec W_n$  separating  two vertices  from $A_n$, and so by \ref{item_kconnect} there again exist at least $k$ edges oriented from $Z$ to $Y$, and also  from $Y$ to $Z$ in $\vec W_n$, and hence in $\vec H$. Together, it follows from Menger's theorem that $\vec H$ is indeed $k$-arc-connected between  vertices of $A_{n+1}$.

Finally, we now lift this orientation of $\vec H$ to an orientation $\vec W_{n+1}$ of the immersion $W_{n+1}$ so that $\vec W_{n+1}$ satisfies \ref{item_exhaust}--\ref{item_kconnect}. Indeed, for each oriented edge in $\vec H$, we simply orient the corresponding path in the immersion $W_{n+1}$ accordingly. Then $\vec W_n \subseteq \vec W_{n+1}$ as directed graphs, and \ref{item_exhaust} holds by construction. To see that property~\ref{item_Eulerian} holds, note that the edges incident with a vertex $v$ in $V(W_{n+1}) \setminus (A_{n+1} \cup X)$ belong to a collection of edge-disjoint, forwards oriented paths containing $v$ in their interior, and hence have equal in- and out-degree. And if for a component $B$ there is a vertex $x$ in $B \cap X$ of degree $3$ in $H$, then since $H$ is $3$-edge-connected and $\vec H$ is well-balanced, it follows that there is at least one ingoing and one outgoing edge at $x$ in $H$, and so $x$ has a difference of $1$ between in- and out-degree in $\vec{W}_{n+1}$. Finally, property \ref{item_kconnect} follows at once from the fact that $\vec W_{n+1}$ is an immersion of the graph $\vec H$ which was $k$-arc-connected between vertices in $A_{n+1}$. 
\end{proof}

\bibliographystyle{plain}
\bibliography{reference}

\begin{bibdiv}
\begin{biblist}

\bib{aharoni1989infinite}{article}{
      author={Aharoni, Ron},
      author={Thomassen, Carsten},
       title={Infinite, highly connected digraphs with no two arc-disjoint
  spanning trees},
        date={1989},
     journal={Journal of graph theory},
      volume={13},
      number={1},
       pages={71\ndash 74},
}

\bib{assem2022analysis}{article}{
      author={Assem, Amena},
       title={Analysis of the lifting graph},
        date={2022},
     journal={arXiv preprint arXiv:2212.03347},
}

\bib{assem2023towards}{article}{
      author={Assem, Amena},
       title={Towards {N}ash-{W}illiams orientation conjecture for infinite
  graphs},
        date={2023},
     journal={arXiv preprint arXiv:2306.10631},
}

\bib{barat2010his}{article}{
      author={Bar{\'a}t, J{\'a}nos},
      author={Kriesell, Matthias},
       title={What is on his mind?},
        date={2010},
     journal={Discrete Mathematics},
      volume={310},
      number={20},
       pages={2573\ndash 2583},
}

\bib{bruhn2007end}{article}{
      author={Bruhn, Henning},
      author={Stein, Maya},
       title={On end degrees and infinite cycles in locally finite graphs},
        date={2007},
     journal={Combinatorica},
      volume={27},
       pages={269\ndash 291},
}

\bib{Bible}{book}{
      author={Diestel, Reinhard},
       title={{Graph Theory}},
     edition={5},
   publisher={Springer},
        date={2015},
}

\bib{egyed1941ueber}{article}{
      author={Egyed, L.},
       title={Ueber die wohlgerichteten unendlichen {G}raphen},
        date={1941},
     journal={Math. phys. Lapok},
      volume={48},
       pages={505\ndash 509},
}

\bib{frank}{article}{
      author={Frank, A.},
       title={On a theorem of {M}ader},
        date={1992},
     journal={Discrete Mathematics},
      volume={101},
       pages={49\ndash 57},
}

\bib{frank2011connections}{book}{
      author={Frank, Andr{\'a}s},
       title={Connections in combinatorial optimization},
   publisher={Oxford University Press Oxford},
        date={2011},
      volume={38},
}

\bib{jannasch2019topologische}{article}{
      author={Jannasch, Lena},
       title={Eine topologische {E}rweiterung des {O}rientierungs-{T}heorems
  von {N}ash-{W}illiams auf lokal endliche {G}raphen},
        date={2019},
     journal={BSc thesis, University of Hamburg},
}

\bib{doi:10.1137/S0895480199364483}{article}{
      author={Jord\'{a}n, Tibor},
       title={Constrained edge-splitting problems},
        date={2003},
     journal={SIAM Journal on Discrete Mathematics},
      volume={17},
      number={1},
       pages={88\ndash 102},
}

\bib{kiraly2006simultaneous}{article}{
      author={Kir{\'a}ly, Zolt{\'a}n},
      author={Szigeti, Zolt{\'a}n},
       title={Simultaneous well-balanced orientations of graphs},
        date={2006},
     journal={Journal of Combinatorial Theory, Series B},
      volume={96},
      number={5},
       pages={684\ndash 692},
}

\bib{lovasz1976some}{article}{
      author={Lov{\'a}sz, L{\'a}szl{\'o}},
       title={On some connectivity properties of eulerian graphs},
        date={1976},
     journal={Acta Mathematica Hungarica},
      volume={28},
      number={1-2},
       pages={129\ndash 138},
}

\bib{mader1978reduction}{incollection}{
      author={Mader, Wolfgang},
       title={A reduction method for edge-connectivity in graphs},
        date={1978},
   booktitle={Annals of {D}iscrete {M}athematics},
      volume={3},
   publisher={Elsevier},
       pages={145\ndash 164},
}

\bib{nash1960orientations}{article}{
      author={Nash-Williams, C.{St}.J.A.},
       title={On orientations, connectivity and odd-vertex-pairings in finite
  graphs},
        date={1960},
     journal={Canadian Journal of Mathematics},
      volume={12},
       pages={555\ndash 567},
}

\bib{nash1969well}{article}{
      author={Nash-Williams, C.{St}.J.A.},
       title={Well-balanced orientations of finite graphs and unobtrusive
  odd-vertex-pairings},
        date={1969},
     journal={Recent Progress in Combinatorics},
       pages={133\ndash 149},
}

\bib{ok2016liftings}{article}{
      author={Ok, Seongmin},
      author={Richter, Bruce},
      author={Thomassen, Carsten},
       title={Liftings in finite graphs and linkages in infinite graphs with
  prescribed edge-connectivity},
        date={2016},
     journal={Graphs and Combinatorics},
      volume={32},
      number={6},
       pages={2575\ndash 2589},
}

\bib{thomassen1989configurations}{inproceedings}{
      author={Thomassen, Carsten},
       title={Configurations in graphs of large minimum degree, connectivity,
  or chromatic number},
        date={1989},
   booktitle={Proceedings of the third international conference on
  combinatorial mathematics},
       pages={402\ndash 412},
}

\bib{thomassen2016orientations}{article}{
      author={Thomassen, Carsten},
       title={Orientations of infinite graphs with prescribed
  edge-connectivity},
        date={2016},
     journal={Combinatorica},
      volume={36},
      number={5},
       pages={601\ndash 621},
}

\end{biblist}
\end{bibdiv}
\end{document}